# ZERO TEMPERATURE LIMIT FOR INTERACTING BROWNIAN PARTICLES. I. MOTION OF A SINGLE BODY[1]

By Tadahisa Funaki

*University of Tokyo*

We consider a system of interacting Brownian particles in $\mathbb{R}^d$ with a pairwise potential, which is radially symmetric, of finite range and attains a unique minimum when the distance of two particles becomes $a > 0$. The asymptotic behavior of the system is studied under the zero temperature limit from both microscopic and macroscopic aspects. If the system is rigidly crystallized, namely if the particles are rigidly arranged in an equal distance $a$, the crystallization is kept under the evolution in macroscopic time scale. Then, assuming that the crystal has a definite limit shape under a macroscopic spatial scaling, the translational and rotational motions of such shape are characterized.

**1. Introduction.** This paper is concerned with a certain scaling limit for a finite, but very large, system of interacting Brownian particles in $\mathbb{R}^d$. The positions of $N$ particles at time $t$, which are denoted by $\mathbf{x}(t) = (x_i(t))_{i=1}^N \in (\mathbb{R}^d)^N$, evolve according to the stochastic differential equation (SDE)

$$(1.1) \qquad dx_i(t) = -\frac{\beta}{2} \nabla_{x_i} H(\mathbf{x}(t))\, dt + dw_i(t), \qquad 1 \le i \le N,$$

where $(w_i(t))_{i=1}^N$ is a family of independent $d$-dimensional standard Brownian motions. The Hamiltonian $H(\mathbf{x})$ of the configuration $\mathbf{x} = (x_i)_{i=1}^N \in (\mathbb{R}^d)^N$ is defined as a sum of pairwise interactions between particles:

$$(1.2) \qquad H(\mathbf{x}) = \sum_{1 \le i < j \le N} U(x_i - x_j).$$

The potential $U = U(x), x \in \mathbb{R}^d$, is radially symmetric, smooth, of finite range and has a unique nondegenerate minimum at $|x| = a > 0$; see Assumption I stated in Section 2 for details. The gradient $\nabla_{x_i} H(\mathbf{x}) \equiv \sum_{j \ne i} \nabla U(x_i -$

Received November 2002; revised March 2003.
[1]Supported in part by JSPS Grants (B)(1)14340029 and 13874015.

$x_j) \in \mathbb{R}^d$ is taken in the variable $x_i$. The parameter $\beta > 0$ represents the inverse temperature of the system.

The basic scaling parameter $\varepsilon > 0$ is the ratio of the microscopic spatial unit length to the macroscopic one. The configuration $\mathbf{x} = (x_i)_{i=1}^N$ is a microscopic object and its macroscopic correspondence is given by $(\varepsilon x_i)_{i=1}^N$ under the spatial scaling $x \mapsto \varepsilon x$. The goal of this paper is to investigate the asymptotic behavior as $\varepsilon \downarrow 0$ of the system defined by (1.1) with $N$, $\beta$ and $t$ suitably scaled depending on $\varepsilon$, especially when the temperature $\beta^{-1}$ of the system converges sufficiently fast to 0. The time change from the microscopic to the macroscopic levels will be introduced for $\mathbf{x}(t)$ by

$$(1.3) \qquad \mathbf{x}^{(\varepsilon)}(t) = \mathbf{x}(\varepsilon^{-\kappa} t), \qquad t \geq 0, \qquad \kappa = d + 2.$$

We say a *rigid crystal* is formed at the microscopic level, if the particles are arranged in an equal distance $a$ and the total energy $H(\mathbf{x})$ increases under any deformation for such arrangement except isometric transformations; see Section 2. As $\beta^{-1} \downarrow 0$, that is, under the zero temperature limit, the system is expected to be frozen and rigidly crystallized.

The results of this paper are twofold and will be formulated at both the microscopic and macroscopic levels in space. The result at the microscopic level can be roughly stated as follows. If the initial configuration $\mathbf{x}(0)$ of the system is nearly a rigid crystal, so is for $\mathbf{x}^{(\varepsilon)}(t)$ asymptotically with probability one as $\varepsilon \downarrow 0$ if the temperature of the system decreases to 0 sufficiently fast; see Theorem 3.4.

The motion of the crystal at the macroscopic level is observed under the spatial scaling $x \mapsto \varepsilon x$. Assuming that the particles' number $N \equiv N(\varepsilon)$ behaves as $\bar{\rho} \varepsilon^{-d}$ with a fixed $\bar{\rho} > 0$ and the crystal has a limit density $\rho(y), y \in \mathbb{R}^d$ as $\varepsilon \downarrow 0$ under the spatial scaling at time $t = 0$, we shall prove that $\mathbf{x}^{(\varepsilon)}(t)$ also has a limit density $\rho_t(y)$ for $t > 0$, which actually coincides with the initial density being isometrically transformed so that $\rho_t(y) = \rho(\varphi_{\theta(t),\eta(t)}^{-1}(y))$ for some $\theta(t)$ and $\eta(t)$. Here, $\varphi_{\theta,\eta}$ denotes an isometry on $\mathbb{R}^d$ defined by $\varphi_{\theta,\eta}(y) = \theta y + \eta, y \in \mathbb{R}^d$ for $\theta = (\theta^{\alpha\beta})_{\alpha,\beta=1}^d \in SO(d), \eta = (\eta^\alpha)_{\alpha=1}^d \in \mathbb{R}^d$; $SO(d)$ stands for the $d$-dimensional special orthogonal group. In other words, the macroscopic limit of $\mathbf{x}^{(\varepsilon)}(t)$ is a rigid body with density $\rho(y)$, which is congruent to the initial body. The translational and rotational motions $(\eta(t), \theta(t))$ of the limit body are random and mutually independent. They are characterized as follows:

$$\eta(t) = (d\text{-dimensional Brownian motion})/\sqrt{\bar{\rho}},$$

while $\theta(t)$ is a Brownian motion on $SO(d)$ which is a solution of an SDE of Stratonovich's type

$$d\theta(t) = \theta(t) \circ dm(t).$$

Here $m(t) = (m^{\alpha\beta}(t))_{\alpha,\beta=1}^d$ satisfies $m^{\alpha\beta}(t) = -m^{\beta\alpha}(t)$ and the upper half components $\{m^{\alpha\beta}(t); \alpha < \beta\}$ of the matrix $m(t)$ are mutually independent such that

$$m^{\alpha\beta}(t) = (\text{one-dimensional Brownian motion})/\sqrt{\bar{q}^\alpha + \bar{q}^\beta}$$

with $\bar{q}^\alpha = \int_{\mathbb{R}^d}(y^\alpha)^2 \rho(y)\,dy$, when the coordinate $y = (y^\alpha)_{\alpha=1}^d$ of $\mathbb{R}^d$ is chosen in such a manner that $\int_{\mathbb{R}^d} y\rho(y)\,dy = 0$ and $(\int_{\mathbb{R}^d} y^\alpha y^\beta \rho(y)\,dy)_{\alpha,\beta=1}^d$ is a diagonal matrix; see Theorem 4.3 and Corollary 4.4. The constants $\bar{\rho}$ and $\bar{q}^\alpha + \bar{q}^\beta$ represent the total mass and moments of inertia [13] of the rigid body with density $\rho(y)$, respectively; note that $\bar{\rho} = \int_{\mathbb{R}^d} \rho(y)\,dy$ holds.

In Section 2, the notion of rigid and infinitesimally rigid crystals is introduced together with several examples. Main results are formulated and proved in Sections 3 and 4. The reason that the time scaling (1.3) is the right one is easily observed for the (macroscopic) translational motion [see the identity (4.15) in the proof of Theorem 4.3], although it may not be obvious for rotation. At the microscopic level, the crystal translates much faster than it rotates. Section 5 is devoted to the proof of technical estimates which are needed for the proof of Theorem 4.3. Section 6 contains concluding remarks.

One of the motivations of this paper comes from the theory of interfaces which appears under the phase transitions. The macroscopic body we have introduced can be regarded as a kind of Wulff shape (see, e.g., [3]) at temperature zero. The static theory for the Wulff shape is recently well developed. This paper attempts to analyze the motion of the Wulff shape by proposing a simple model. Indeed, at temperature zero, at least two pure phases arise in our model. One is the high density region where particles are arranged in an equal distance $a$ and the other is the empty region where there are no particles. In this respect, the body with density $\rho(y)$ is a mixture of high density and empty regions observed macroscopically. It would be more natural and desirable to study the model with temperature being sufficiently small but fixed under the scaling. However, this problem turns out to be quite hard. Actually, to solve such a problem, we need to have information on the phase transition for the Gibbs measures corresponding to our dynamics with infinitely many particles, but it is not well known. Lang [14] considered a system of ordinary differential equation (1.1) dropping Brownian motions with $\beta = 1$ and $N = \infty$. Such a system is obtained from the SDE in the zero temperature limit $\beta \to \infty$ under a time change $t \mapsto \beta^{-1}t$. In one dimension and for strictly convex potential $U$ having a hard core, ergodic properties of the dynamics were studied and equilibrium states (called rigid states) were characterized. Related problems were discussed for the stochastic partial differential equations in [4] and [5].

This paper deals with the motion of a single body (or single crystal in microscopic aspect). The coagulations of several bodies are discussed in [6] in one dimension. The study of coagulations in higher dimensions is out of reach at present.

**2. Rigid configurations.** We introduce the notion of rigidity and infinitesimal rigidity for configurations of particles in $\mathbb{R}^d$ and expose several examples of infinitesimally rigid configurations. The number $N$ of particles is fixed throughout this section.

2.1. *Hamiltonian and rigidity of configurations.* The Hamiltonian $H(\mathbf{x})$ of $\mathbf{x} \in (\mathbb{R}^d)^N$ is introduced by the formula (1.2). The potential $U$ is assumed to be radially symmetric; that is, $U(x) = U(|x|)$, $x \in \mathbb{R}^d$ for $U = U(r), r \geq 0$, and the function $U(r)$ satisfies the following conditions:

ASSUMPTION I. (i) (Smoothness, finite range). $U \in C_0^3(\mathbb{R})$, where $U(-r) := U(r)$.
(ii) There exists a unique $a > 0$ such that $U(a) = \min_{r \geq 0} U(r)$ and $\check{c} := U''(a) > 0$.

Condition (ii) means that the energy for two particles takes minimal value when the distance between these particles becomes $a$. The range of $U$ is defined by $b := \inf\{r > 0; U(s) = 0 \text{ for every } s > r\}$. Let $\mathbf{z} = (z_i)_{i=1}^N \in (\mathbb{R}^d)^N$ be a configuration satisfying

$$|z_i - z_j| = a \quad \text{or} \quad |z_i - z_j| > b \tag{2.1}$$

for every $1 \leq i \neq j \leq N$. An additional condition on $b$ is necessary for such $\mathbf{z}$ to exist; for example, see condition (2.6) for configurations on a triangular lattice. The configuration $\mathbf{z}$ is a critical point of the Hamiltonian $H$. This physically means that $\mathbf{z}$ is a microscopically crystallized or frozen configuration of atoms at temperature zero. Its rotated and shifted configuration $\varphi_{\theta,\eta}(\mathbf{z}) := (\varphi_{\theta,\eta}(z_i))_{i=1}^N \equiv (\theta z_i + \eta)_{i=1}^N$ is obviously a critical point of $H$ again for every $\theta \in SO(d)$ and $\eta \in \mathbb{R}^d$. We shall write

$$\mathcal{M} = \{\varphi_{\theta,\eta}(\mathbf{z}); \theta \in SO(d), \eta \in \mathbb{R}^d\} \subset (\mathbb{R}^d)^N, \tag{2.2}$$

and its tubular neighborhood

$$\mathcal{M}_2(\delta) = \{\mathbf{x} \in (\mathbb{R}^d)^N; \text{dist}(\mathbf{x}, \mathcal{M}) \leq \delta\}, \qquad \delta > 0, \tag{2.3}$$

where the distance is defined under the Euclidean norm in $(\mathbb{R}^d)^N$: $\text{dist}(\mathbf{x}, \mathcal{M}) = \inf_{\mathbf{y} \in \mathcal{M}} \|\mathbf{x} - \mathbf{y}\|_2$ and $\|\mathbf{x} - \mathbf{y}\|_2 = (\sum_{i=1}^N |x_i - y_i|^2)^{1/2}$. The configuration $\mathbf{z}$ satisfying (2.1) will be called a *crystal*.

We say the crystal $\mathbf{z}$ is *rigid* if $H(\mathbf{x}) > H(\mathbf{z})$ holds for every $\mathbf{x} \in \mathcal{M}_2(\delta) \setminus \mathcal{M}$ and some $\delta > 0$. The rigidity means that $\mathbf{z}$ has no internal degree of freedom except for the isometry. For example, in two dimension, the three vertices of an equilateral triangle form a rigid crystal, but the four vertices of a square do not. The rigid crystal is a local minimum of $H$ by definition, but not necessarily a global one.

2.2. *Orthogonal decomposition of* $\mathbf{x} \in \mathcal{M}_2(\delta)$. In order to study the rigidity of a crystal $\mathbf{z}$, we introduce a decomposition of $\mathbf{x} \in \mathcal{M}_2(\delta)$. Let $\mathcal{H}_{\mathbf{z}} = \{X\mathbf{z} + h; X \in \mathfrak{so}(d), h \in \mathbb{R}^d\} \subset (\mathbb{R}^d)^N$ be the tangent space to $\mathcal{M}$ at $\mathbf{z}$, where $X\mathbf{z} + h := (Xz_i + h)_{i=1}^N$ and $\mathfrak{so}(d) = \{X \in M(d); X + {}^t X = 0\}$ is the Lie algebra of $SO(d)$. The set $M(d)$ stands for the family of all $d \times d$ real matrices. The orthogonal subspace to $\mathcal{H}_{\mathbf{z}}$ in $(\mathbb{R}^d)^N$ under the inner product $(\mathbf{h}, \mathbf{h}') := \sum_{i=1}^N (h_i, h_i')$ for $\mathbf{h} = (h_i)_{i=1}^N$ and $\mathbf{h}' = (h_i')_{i=1}^N$ is denoted by $\mathcal{H}_{\mathbf{z}}^\perp$. The corresponding norm of $\mathbf{h}$ is $\|\mathbf{h}\|_2 = (\sum_{i=1}^N |h_i|^2)^{1/2}$ defined previously.

For every $\mathbf{x} \in \mathcal{M}_2(\delta)$, we denote by $\mathbf{z}(\mathbf{x}) := \mathbf{y} \in \mathcal{M}$ the minimizer of $\|\mathbf{x} - \mathbf{y}\|_2$ in $\mathbf{y} \in \mathcal{M}$. If $\delta > 0$ is sufficiently small, $\mathbf{z}(\mathbf{x})$ is uniquely determined and $\mathbf{x} \in \mathcal{M}_2(\delta)$ admits a decomposition:

$$(2.4) \qquad \mathbf{x} = \mathbf{z}(\mathbf{x}) + \mathbf{h}(\mathbf{x}), \qquad \mathbf{z}(\mathbf{x}) \in \mathcal{M}, \ \mathbf{h}(\mathbf{x}) \in \mathcal{H}_{\mathbf{z}(\mathbf{x})}^\perp.$$

In fact, by the definition of $\mathbf{z}(\mathbf{x})$, we have

$$\left. \frac{d}{du} \|\mathbf{x} - \varphi_{e^{uX}, uh}(\mathbf{z}(\mathbf{x}))\|_2^2 \right|_{u=0} = 0$$

for every $X \in \mathfrak{so}(d)$ and $h \in \mathbb{R}^d$, and this implies $\mathbf{h}(\mathbf{x}) := \mathbf{x} - \mathbf{z}(\mathbf{x}) \in \mathcal{H}_{\mathbf{z}(\mathbf{x})}^\perp$.

2.3. *Hessian of $H$ on $\mathcal{M}$.* Let

$$\operatorname{Hess} H(\mathbf{x}) = \left( \frac{\partial^2 H}{\partial x_i^\alpha \partial x_j^\beta}(\mathbf{x}) \right)_{1 \leq \alpha, \beta \leq d, \ 1 \leq i, j \leq N} \in M(dN)$$

be the Hessian of $H$ and define a quadratic form in $\mathbf{h} = (h_i)_{i=1}^N \in (\mathbb{R}^d)^N$ by

$$(\mathbf{h}, \operatorname{Hess} H(\mathbf{x})\mathbf{h}) = \sum_{i,j=1}^N \sum_{\alpha, \beta=1}^d \frac{\partial^2 H}{\partial x_i^\alpha \partial x_j^\beta}(\mathbf{x}) h_i^\alpha h_j^\beta.$$

Then a direct calculation yields the following at $\mathbf{x} = \mathbf{z}$.

LEMMA 2.1.

$$\mathcal{E}_1(\mathbf{h}) \equiv \mathcal{E}_{1,\mathbf{z}}(\mathbf{h}) := (\mathbf{h}, \operatorname{Hess} H(\mathbf{z})\mathbf{h}) = \frac{\check{c}}{a^2} \sum_{\langle i,j \rangle} (h_i - h_j, z_i - z_j)^2 \geq 0,$$

*where the sum $\langle i, j \rangle$ is taken over all pairs $\{i, j\}$ satisfying $|z_i - z_j| = a$. We call such pairs neighboring.*

This lemma immediately shows that the Hessian of $H$ degenerates for $\mathbf{h} \in \mathcal{H}_{\mathbf{z}}$. The degeneracy for $(h_i = h)_{i=1}^N$ comes from the invariance of $H$ under the translation, while that for $(h_i = Xz_i)_{i=1}^N$ comes from its invariance under the rotation. The rigidity of $\mathbf{z}$ follows from the nondegeneracy of the Hessian of $H$ for $\mathbf{h} \in \mathcal{H}_{\mathbf{z}}^\perp$. This leads us to introduce the following notion.

DEFINITION 2.1. We call the crystal $\mathbf{z}$ infinitesimally rigid if the quadratic form Hess $H(\mathbf{z})$ restricted on the subspace $\mathcal{H}_{\mathbf{z}}^{\perp}$ is (strictly) positive definite: $\mathcal{E}_1(\mathbf{h}) = 0 \Longleftrightarrow \mathbf{h} \in \mathcal{H}_{\mathbf{z}}$.

Since $H(\mathbf{x}) = H(\mathbf{z}) + \mathcal{E}_1(\mathbf{h}(\mathbf{x}))/2 + o(\|\mathbf{h}(\mathbf{x})\|_2^2)$ as $\|\mathbf{h}(\mathbf{x})\|_2 \to 0$ under the decomposition (2.4) [see (3.1)], we easily see that the infinitesimal rigidity implies the rigidity.

REMARK 2.1. The study of rigidity and infinitesimal rigidity for bar and joint frameworks has a long history; see [1], [2] and [17]. The length of bars is always $a$ in our case, but such an assumption is unnecessary in a general theory. According to ([1], page 281), $\mathbf{z}$ is called rigid in $\mathbb{R}^d$ if, for every $\mathbf{x}$ sufficiently close to $\mathbf{z}$ satisfying $|x_i - x_j| = a$ for every neighboring pair $\langle i, j \rangle$, there exists an isometry $\varphi$ of $\mathbb{R}^d$ such that $x_i = \varphi(z_i)$ holds for every $i$. In [2] (page 173) $\mathbf{z}$ is called infinitesimally rigid in $\mathbb{R}^d$ if $T_{\mathbf{z}} = \ker df_G(\mathbf{z})$ holds, where $T_{\mathbf{z}} = \mathcal{H}_{\mathbf{z}}$ and $\ker df_G(\mathbf{z}) = \{\mathbf{h} \in (\mathbb{R}^d)^N; \mathcal{E}_1(\mathbf{h}) = 0\}$, respectively, in our terminology. We therefore see that the definitions of rigidity and infinitesimal rigidity employed by these papers coincide with ours. Three points $\{p_1, p_2, p_3\}$ in $\mathbb{R}^2$ sitting on a line are rigid but not infinitesimally rigid in $\mathbb{R}^2$ when the distances between any two points are specified; see [12], noting that the rigidity in that paper means the infinitesimal one. This example is not for a crystal, but exhibits the difference in two notions.

The rigidity of $\mathbf{z}$ implies the connectedness of the set $\mathbf{z}$ under the neighboring relation $\langle i, j \rangle$ and therefore we have, under the infinitesimal rigidity, the spectral gap for the quadratic form $\mathcal{E}_1(\mathbf{h})$ in the following sense:

$$\lambda^{(1)}(\mathbf{z}) = \inf\left\{\frac{\mathcal{E}_1(\mathbf{h})}{\|\nabla \mathbf{h}\|_2^2}; \mathbf{h} \in \mathcal{H}_{\mathbf{z}}^{\perp}, \|\nabla \mathbf{h}\|_2 \neq 0\right\} > 0,$$

where

$$\|\nabla \mathbf{h}\|_2^2 = \sum_{\langle i,j \rangle} |h_i - h_j|^2.$$

This can be rewritten as

(2.5) $\qquad \lambda^{(1)}(\mathbf{z})\|\nabla \mathbf{h}\|_2^2 \leq \mathcal{E}_1(\mathbf{h}) \leq \check{c}\|\nabla \mathbf{h}\|_2^2, \qquad \mathbf{h} \in \mathcal{H}_{\mathbf{z}}^{\perp}.$

Note that the second inequality is obvious.

REMARK 2.2. In one dimension, a chain $\mathbf{z} = (z_i = ai)_{i=1}^N$ arranged in an equal distance $a$ is a rigid crystal and $\lambda^{(1)}(\mathbf{z}) = \check{c}$ since $\mathcal{E}_1(\mathbf{h}) = \check{c}\|\nabla \mathbf{h}\|_2^2$.

2.4. *Examples of infinitesimally rigid crystals.* We prepare two lemmas before constructing several examples of infinitesimally rigid crystals. We say a set $C = \{x_i\}_{i=0}^n \subset \mathbb{R}^d, n \leq d$, is an $n$-dimensional cell in $\mathbb{R}^d$ if the dimension of the affine hull of $C$ is $n$ and $|x_i - x_j| = a$ for every $0 \leq i < j \leq n$.

LEMMA 2.2. *Let $\{e_i\}_{i=1}^d$ be a basis of $\mathbb{R}^d$ and set $e_0 = 0$. If vectors $\{h_i\}_{i=0}^d \subset \mathbb{R}^d$ satisfy $(h_i - h_j, e_i - e_j) = 0$ for every $0 \leq i < j \leq d$, then there exists a unique $X \in \mathfrak{so}(d)$ such that $h_i = Xe_i + h_0, 1 \leq i \leq d$. In particular, a $d$-dimensional cell is infinitesimally rigid.*

PROOF. The uniqueness of $X$ is obvious, since $\{e_i\}_{i=1}^d$ forms a basis of $\mathbb{R}^d$. To show the existence of $X$, we may assume that $(e_1, \ldots, e_d)$ is an upper triangular matrix. In fact, by Schmidt's orthogonalization, one can find $P \in O(d)$ and upper triangular matrix $(\tilde{e}_1, \ldots, \tilde{e}_d)$ such that $(e_1, \ldots, e_d) = P(\tilde{e}_1, \ldots, \tilde{e}_d)$. If the conclusion holds for $(\tilde{e}_1, \ldots, \tilde{e}_d)$, there exists $\tilde{X} \in \mathfrak{so}(d)$ such that $P^{-1}x_i = \tilde{X}\tilde{e}_i, 1 \leq i \leq d$. Taking $X = P\tilde{X}P^{-1} \in \mathfrak{so}(d)$, the conclusion is shown also for $\{e_i\}_{i=1}^d$.

Now we assume $(e_1, \ldots, e_d)$ is an upper triangular matrix and use an induction in $d$ to construct $X$. We may further assume $h_0 = 0$ by replacing $h_i$ with $h_i - h_0$. Since $e_i$ has a form

$$e_i = \begin{pmatrix} e_i' \\ 0 \end{pmatrix}$$

with $e_i' \in \mathbb{R}^{d-1}$ for $1 \leq i \leq d-1$, writing

$$h_i = \begin{pmatrix} h_i' \\ h_i'' \end{pmatrix}$$

with $h_i' \in \mathbb{R}^{d-1}$ and $h_i'' \in \mathbb{R}$, we have

$$0 = (h_i - h_j, e_i - e_j) = (h_i' - h_j', e_i' - e_j'), \qquad 1 \leq i < j \leq d-1,$$

and therefore $h_i' = X'e_i'$ holds for some $X' \in \mathfrak{so}(d-1)$ and every $1 \leq i \leq d-1$ by the assumption of the induction. Now, writing

$$e_d = \begin{pmatrix} e_d' \\ e_d'' \end{pmatrix}$$

with $e_d' \in \mathbb{R}^{d-1}$ and $e_d'' \in \mathbb{R} \setminus \{0\}$, define

$$X = \begin{pmatrix} X' & y \\ -{}^t y & 0 \end{pmatrix} \in \mathfrak{so}(d) \quad \text{with} \quad y = \frac{1}{e_d''}(h_d' - X'e_d') \in \mathbb{R}^{d-1}.$$

Then, we see that $h_d = Xe_d$ holds. We need to prove that $h_i = Xe_i$ holds also for $1 \leq i \leq d-1$. To this end, it is enough to show that $-(y, e_i') = h_i''$.

However, since $(h_i - h_d, e_i - e_d) = (h_i, e_i) = (h_d, e_d) = 0$, we have $(h_i, e_d) + (h_d, e_i) = 0$. This implies $-(y, e'_i) = h''_i$, since $(h_d, e_i) = (h'_d, e'_i)$ and

$$(h_i, e_d) = (h'_i, e'_d) + h''_i e''_d = (X'e'_i, e'_d) + h''_i e''_d$$
$$= -(e'_i, X'e'_d) + h''_i e''_d = -(e'_i, h'_d) + e''_d(e'_i, y) + h''_i e''_d.$$

Therefore, the conclusion is shown in $d$ dimension if it is true in $d-1$ dimension. The procedure of the induction is complete once we can show the conclusion when $d = 2$. However, this is already essentially done in the above argument.

The infinitesimal rigidity of $d$-dimensional cell $C = \{x_i\}_{i=0}^d$ is immediate by taking $e_i = x_i - x_0$ for $1 \leq i \leq d$. □

LEMMA 2.3. *Let two infinitesimally rigid crystals $\mathbf{z}^{(1)}$ and $\mathbf{z}^{(2)}$ be given and assume that the dimension of the affine hull of $\mathbf{z}^{(1)} \cap \mathbf{z}^{(2)}$ is at least $d-1$. Then, the joined configuration $\mathbf{z}^{(1)} \cup \mathbf{z}^{(2)}$ is infinitesimally rigid.*

PROOF. Let us denote $\mathbf{z}^{(1)} \cap \mathbf{z}^{(2)} = \{z_i^{(0)}\}_{i=1}^{N_0}$. Then, the conclusion follows if one can show that $X z_i^{(0)} + h = X' z_i^{(0)} + h'$ for every $1 \leq i \leq N_0$ implies $X = X'$ and $h = h'$, where $X, X' \in \mathfrak{so}(d)$ and $h, h' \in \mathbb{R}^d$. However, from the assumption, one can find at least $d-1$ linearly independent vectors $\{e_k\}_{k=1}^{d-1}$ from $\{z_i^{(0)} - z_j^{(0)}\}_{1 \leq i < j \leq N_0}$. The identities $Xe_k = X'e_k$ hold for such vectors. Take $e_d \in \mathbb{R}^d$ in such a manner that $\{e_k\}_{k=1}^d$ forms a basis of $\mathbb{R}^d$. Then, since

$$(Xe_d, e_k) = -(e_d, Xe_k) = -(e_d, X'e_k) = (X'e_d, e_k)$$

for $1 \leq k \leq d-1$ and $(Xe_d, e_d) = 0 = (X'e_d, e_d)$, we see $Xe_d = X'e_d$. We accordingly have $Xe_k = X'e_k$ for every $1 \leq k \leq d$. This proves $X = X'$ and therefore $h = h'$. □

EXAMPLE 2.1. The set obtained by patching together $d$-dimensional cells on their faces is infinitesimally rigid from Lemmas 2.2 and 2.3. More precisely, a finite set $A \subset \mathbb{R}^d$ satisfying the following two conditions is infinitesimally rigid:

1. $A = \bigcup_k C_k$ with finitely many $d$-dimensional cells $C_k$.
2. For any $C_{k_1}$ and $C_{k_n}$ in $A$, there exists a sequence $C_{k_2}, \ldots, C_{k_{n-1}}$ in $A$ such that $C_{k_i} \cap C_{k_{i+1}}$ are $(d-1)$-dimensional cells for $1 \leq i \leq n-1$.

In two dimension, the set as in Figure 1 is infinitesimally rigid. In general,

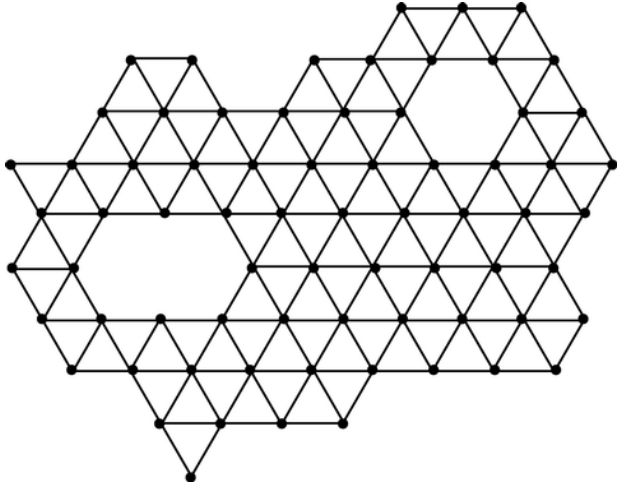

Fig. 1. *Two-dimensional crystal.*

infinitesimally rigid crystals may have defects.

EXAMPLE 2.2. In three dimension, the tetrahedron ( = three-dimensional cell), octahedron and icosahedron are infinitesimally rigid by Cauchy's rigidity theorem or by Alexandrov's rigidity theorem [2]. Note that the faces of the three types of regular polyhedrons listed above are all equilateral triangles. In particular, the set obtained by patching together polyhedrons of these three types as in Example 2.1 is infinitesimally rigid.

EXAMPLE 2.3 (Crystals on triangular lattice). Let $\{e_\alpha \in \mathbb{R}^d\}_{\alpha=1}^d$ be a basis of $\mathbb{R}^d$ such that $(e_\alpha, e_\beta) = (1 + \delta_{\alpha\beta})/2$. In other words, it is a system of unit vectors and arbitrarily chosen two of them are at an angle of $60°$ with each other. Then, a $d$-dimensional triangular lattice is defined as an integer lattice generated by $\{e_\alpha\}_{\alpha=1}^d$: $\Lambda \equiv \Lambda_d = \{\sum_{\alpha=1}^d \xi_\alpha e_\alpha \in \mathbb{R}^d; \xi = (\xi_\alpha)_{\alpha=1}^d \in \mathbb{Z}^d\}$. Note that $\Lambda_1 = \mathbb{Z}$. Set $E = \{e \in \Lambda; |e| = 1\}$, and then it is easy to see that $e = \sum_{\alpha=1}^d \xi_\alpha e_\alpha \in E$ if and only if $\xi$ has the form $\xi_\alpha = \pm\delta_{\alpha_0\alpha}$ for some $\alpha_0$ or $\xi_\alpha = \delta_{\alpha_0\alpha} - \delta_{\beta_0\alpha}$ for some $\alpha_0 \ne \beta_0$. The triangular lattice $\Lambda$ is the set of centers of circles (when $d = 2$) or balls (when $d = 3$) with radius $1/2$ filled most densely in the space; this assertion was known as the Kepler conjecture in three dimension and solved by Hales [8]. We need an additional assumption:

(2.6) $$b < c(\Lambda_d)a,$$

for a rigid **z** to exist on $a\Lambda$ satisfying (2.1), where $c(\Lambda_d) := \inf_{x \in \Lambda_d \setminus E} |x| = 2$ (when $d = 1$), $\sqrt{3}$ (when $d = 2$) and $\sqrt{2}$ (when $d \ge 3$).

In two dimension, $\Lambda$ can be constructed by patching equilateral triangles, while in three dimension, $\Lambda$ is obtainable by patching tetrahedrons and octahedrons. Therefore, at least in two and three dimensions, $\varepsilon^{-1} D \cap a\Lambda$ is infinitesimally rigid for a bounded domain $D$ in $\mathbb{R}^d$ ($d = 2$ or $3$) having smooth boundary and for small $\varepsilon$ by deleting or adding some points near the boundary in a proper way if necessary.

2.5. *Tubular neighborhoods of $\mathcal{M}$ defined in two other norms.* Let us consider two norms $\|\mathbf{h}\|_\infty$ and $\|\nabla \mathbf{h}\|_\infty$ for $\mathbf{h} = (h_i)_{i=1}^N \in (\mathbb{R}^d)^N$ defined by

$$\|\mathbf{h}\|_\infty = \max_i |h_i| \quad \text{and} \quad \|\nabla \mathbf{h}\|_\infty = \max_{\langle i,j \rangle} |h_i - h_j|,$$

respectively; recall that $\langle i, j \rangle$ refers to neighboring pairs. Then, for every small $c > 0$, tubular neighborhoods $\mathcal{M}_\infty(c)$ and $\mathcal{M}_\infty^\nabla(c)$ of $\mathcal{M}$ can be introduced as

$$(2.7) \qquad \mathcal{M}_\infty(c) = \{\mathbf{x} \in (\mathbb{R}^d)^N; \|\mathbf{h}(\mathbf{x})\|_\infty \leq c\},$$

$$(2.8) \qquad \mathcal{M}_\infty^\nabla(c) = \{\mathbf{x} \in (\mathbb{R}^d)^N; \|\nabla \mathbf{h}(\mathbf{x})\|_\infty \leq c\},$$

respectively, where $\mathbf{h}(\mathbf{x}) \in (\mathbb{R}^d)^N$ is defined by (2.4) for $\mathbf{x} \in \mathcal{M}_2(\delta)$ with sufficiently small $\delta > 0$. For each crystal $\mathbf{z}$, since two norms $\|\mathbf{h}\|_\infty$ and $\|\mathbf{h}\|_2$ are mutually equivalent, one can find $\bar{c} = \bar{c}(\mathbf{z}) > 0$ such that $\mathbf{h}(\mathbf{x})$ is well defined for all $\mathbf{x} \in \mathcal{M}_\infty(\bar{c}(\mathbf{z}))$.

For $1 \leq i \neq j \leq N$, let $\mathfrak{p}(i,j) = \{i = i_0 \sim i_1 \sim \cdots \sim i_n = j\}$ be the shortest path connecting $i$ and $j$, where $i_k \sim i_{k+1}$ means that the pair $\langle i_k, i_{k+1} \rangle$ is neighboring. We call $n =: \sharp\mathfrak{p}(i,j)$ the length of $\mathfrak{p}(i,j)$ and define the radius of $\mathbf{z}$ by

$$R(\mathbf{z}) = \max\{\sharp \mathfrak{p}(i,j); 1 \leq i \neq j \leq N\}.$$

LEMMA 2.4. *For every $\mathbf{x} \in \mathcal{M}_\infty(\bar{c}(\mathbf{z}))$, we have*

$$(2.9) \qquad \|\mathbf{h}(\mathbf{x})\|_\infty \leq R(\mathbf{z}) \|\nabla \mathbf{h}(\mathbf{x})\|_\infty.$$

*In particular, the set $\mathcal{M}_\infty^\nabla(c)$ is well defined for $0 < c \leq \bar{c}(\mathbf{z})/R(\mathbf{z})$.*

PROOF. Since $\sum_{j=1}^N h_j = 0$ for $\mathbf{h} \equiv \mathbf{h}(\mathbf{x}) = (h_i)_{i=1}^N$,

$$|h_i| = \left| h_i - \frac{1}{N} \sum_{j=1}^N h_j \right| \leq \frac{1}{N} \sum_{j=1}^N |h_i - h_j|$$

$$\leq \frac{1}{N} \sum_{j=1}^N \sum_{\langle i_k, i_{k+1}\rangle \in \mathfrak{p}(i,j)} |h_{i_k} - h_{i_{k+1}}| \leq R(\mathbf{z}) \|\nabla \mathbf{h}(\mathbf{x})\|_\infty$$

for every $1 \leq i \leq N$. This shows (2.9) and $\mathcal{M}_\infty^\nabla(c) \subset \mathcal{M}_\infty(R(\mathbf{z})c)$. $\square$

**3. Microscopic shape theorem.** This section establishes the asymptotic behavior as $\varepsilon \downarrow 0$ of $\mathbf{x}^{(\varepsilon)}(t)$, the solution of the SDE (1.1) which is scaled macroscopically in time, when the initial configuration $\mathbf{x}^{(\varepsilon)}(0)$ is nearly infinitesimally rigid and the temperature $\beta^{-1} = \beta(\varepsilon)^{-1}$ of the system decreases to 0 sufficiently fast compared with $\varepsilon$.

3.1. *Behavior of $H$ near $\mathcal{M}$.* Let $\mathbf{z}$ be a crystal, that is, a configuration satisfying the condition (2.1). A configuration $\mathbf{x} = (x_i)_{i=1}^N$ is then decomposed as $\mathbf{x} = \mathbf{z} + \mathbf{h}$ around $\mathbf{z}$ just by setting $\mathbf{h} \equiv (h_i)_{i=1}^N := (x_i - z_i)_{i=1}^N$. We shall write

$$G(\mathbf{x}) = \sum_{i=1}^N |\nabla_{x_i} H(\mathbf{x})|^2 (\equiv \|\nabla H(\mathbf{x})\|_2^2).$$

LEMMA 3.1. *Suppose that $\mathbf{x}$ satisfies $|x_i - x_j| > b$ for nonneighboring pairs $\{i,j\}$ (of $\mathbf{z}$) and $|x_i - x_j| \geq a_0$ for neighboring pairs $\langle i,j \rangle$ with some $a_0 \in (0, a)$. Then, there exists $C > 0$ (independent of $\mathbf{x}, \mathbf{z}$ and $N$) such that*

$$(3.1) \qquad |H(\mathbf{x}) - \{H(\mathbf{z}) + \tfrac{1}{2}\mathcal{E}_1(\mathbf{h})\}| \leq C \sum_{\langle i,j \rangle} |h_i - h_j|^3,$$

$$(3.2) \qquad |G(\mathbf{x}) - \mathcal{E}_2(\mathbf{h})| \leq C \sum_{\langle i,j \rangle} |h_i - h_j|^3,$$

*where*

$$\mathcal{E}_2(\mathbf{h}) := \frac{\check{c}^2}{a^4} \sum_{i=1}^N \left| \sum_{j:\langle i,j\rangle} (h_i - h_j, z_i - z_j)(z_i - z_j) \right|^2 = \frac{1}{4} \sum_{i=1}^N |\nabla_{h_i} \mathcal{E}_1(\mathbf{h})|^2.$$

PROOF. These two estimates are shown by Taylor's theorem applied for $H(\mathbf{x})$ and $G(\mathbf{x})$ in the variables $\{h_i - h_j\}$ with neighboring pairs $\langle i,j \rangle$. Since the computations are easy, the details are omitted. □

Now let us assume that $\mathbf{z}$ is infinitesimally rigid.

LEMMA 3.2. *There exist $C > 0$ and $\lambda^{(2)}(\mathbf{z}) > 0$ such that*

$$C^{-1}\mathcal{E}_2(\mathbf{h}) \leq \mathcal{E}_1(\mathbf{h}) \leq \{\lambda^{(2)}(\mathbf{z})\}^{-1}\mathcal{E}_2(\mathbf{h}), \qquad \mathbf{h} \in \mathcal{H}_{\mathbf{z}}^\perp.$$

PROOF. The first inequality is obvious. To show the second, note that the quadratic form $\mathcal{E}_1(\mathbf{h})$ is expressed as $\mathcal{E}_1(\mathbf{h}) = (A\mathbf{h}, \mathbf{h})$ with a symmetric matrix $A \in M(dN)$ and $\mathcal{H}_{\mathbf{z}}$ is the eigenspace of $A$ corresponding to the eigenvalue 0. Since we have

$$\mathcal{E}_2(\mathbf{h}) = \|A\mathbf{h}\|_2^2 = (A^2\mathbf{h}, \mathbf{h}),$$

$\mathcal{E}_2(\mathbf{h}) = 0$ holds if and only if $\mathbf{h} \in \mathcal{H}_{\mathbf{z}}$ and this implies the conclusion. □

REMARK 3.1. The constant $\lambda^{(2)}(\mathbf{z})$ is related to the Poincaré inequality; see Lemma 2.1 of [6] in one dimension. We may assume $0 < \lambda^{(2)}(\mathbf{z}) \leq 1$ for every $\mathbf{z}$.

In the following, we shall normalize the Hamiltonian $H$ as $H(\mathbf{z}) = 0$ by adding a constant [i.e., by considering $H - H(\mathbf{z})$ instead of $H$].

COROLLARY 3.3. *If* $\mathbf{h} = \mathbf{x} - \mathbf{z} \in \mathcal{H}_\mathbf{z}^\perp$ *and* $\|\nabla \mathbf{h}\|_\infty \leq \delta \lambda^{(1)}(\mathbf{z}) \lambda^{(2)}(\mathbf{z})$ *is satisfied for sufficiently small* $\delta > 0$, *we have*

$$C^{-1} \lambda^{(2)}(\mathbf{z}) H(\mathbf{x}) \leq G(\mathbf{x}) \leq C H(\mathbf{x}),$$

*for some* $C > 0$.

PROOF. The estimates (3.1), (2.5) and the assumption on $\|\nabla \mathbf{h}\|_\infty$ show that

$$|H(\mathbf{x}) - \tfrac{1}{2}\mathcal{E}_1(\mathbf{h})| \leq C\|\nabla \mathbf{h}\|_\infty \|\nabla \mathbf{h}\|_2^2 \leq C\delta \mathcal{E}_1(\mathbf{h}),$$

since $\lambda^{(2)}(\mathbf{z}) \leq 1$. Therefore, taking $\delta$ sufficiently small so that $C\delta \leq 1/4$, we have

(3.3) $$\tfrac{1}{4}\mathcal{E}_1(\mathbf{h}) \leq H(\mathbf{x}) \leq \tfrac{3}{4}\mathcal{E}_1(\mathbf{h}).$$

On the other hand, from (3.2), (2.5) and Lemma 3.2,

$$|G(\mathbf{x}) - \mathcal{E}_2(\mathbf{h})| \leq C\|\nabla \mathbf{h}\|_\infty \|\nabla \mathbf{h}\|_2^2 \leq C\delta \mathcal{E}_2(\mathbf{h}),$$

which shows

(3.4) $$\tfrac{3}{4}\mathcal{E}_2(\mathbf{h}) \leq G(\mathbf{x}) \leq \tfrac{5}{4}\mathcal{E}_2(\mathbf{h}).$$

The conclusion follows from (3.3), (3.4) and Lemma 3.2. □

3.2. *Lyapounov argument.* Assume that a sequence $\mathbf{z}^{(\varepsilon)} = (z_i^{(\varepsilon)})_{i=1}^N$, $0 < \varepsilon < 1$, of infinitesimally rigid and centered crystals is given, where "centered" means $\sum_{i=1}^N z_i^{(\varepsilon)} = 0$. The number $N \equiv N(\varepsilon)$ of particles in $\mathbf{z}^{(\varepsilon)}$ may change depending on the scaling parameter $\varepsilon$. Let $\mathbf{x}(t) = (x_i(t))_{i=1}^N \in (\mathbb{R}^d)^N$ be the solution of the SDE (1.1) and introduce the time changed process $\mathbf{x}^{(\varepsilon)}(t)$ of $\mathbf{x}(t)$ by (1.3). The inverse temperature $\beta \equiv \beta(\varepsilon)$ changes with $\varepsilon$ and diverges to $+\infty$ as $\varepsilon \downarrow 0$ sufficiently fast; see the condition (3.8) in Theorem 3.4.

Let $\mathcal{M} \equiv \mathcal{M}^{(\varepsilon)}$ and $\mathcal{M}_\infty^\nabla(c) \equiv \mathcal{M}_\infty^{\nabla,(\varepsilon)}(c)$ be the sets (2.2) and (2.8) determined from $\mathbf{z}^{(\varepsilon)}$ instead of $\mathbf{z}$, respectively. Given $c = c(\varepsilon) \downarrow 0$, consider a sequence of stopping times $\sigma \equiv \sigma^{(\varepsilon)}$ defined by

(3.5) $$\sigma = \inf\{t \geq 0; \mathbf{x}^{(\varepsilon)}(t) \notin \mathcal{M}_\infty^{\nabla,(\varepsilon)}(c(\varepsilon))\}.$$

The main result of this section can now be stated. The proof is based on the Lyapounov type argument.

THEOREM 3.4. *Let $\{c = c(\varepsilon) \downarrow 0\}$ (as $\varepsilon \downarrow 0$) and a sequence of (random) initial data $\{\mathbf{x}^{(\varepsilon)}(0)\}$ be given and satisfy the following conditions:*

$$(3.6) \qquad 0 < c(\varepsilon) \le \delta \lambda^{(1,\varepsilon)} \lambda^{(2,\varepsilon)} \wedge \bar{c}^{(\varepsilon)},$$

$$(3.7) \qquad E[\|\nabla \mathbf{h}(\mathbf{x}^{(\varepsilon)}(0))\|_2^{2p}]^{1/2p} = o(\{\lambda^{(1,\varepsilon)}\}^{1/2} c(\varepsilon)),$$

*as $\varepsilon \downarrow 0$ for some $p > 1$, where $\lambda^{(1,\varepsilon)} = \lambda^{(1)}(\mathbf{z}^{(\varepsilon)})$, $\lambda^{(2,\varepsilon)} = \lambda^{(2)}(\mathbf{z}^{(\varepsilon)})$, $\bar{c}^{(\varepsilon)} = \bar{c}(\mathbf{z}^{(\varepsilon)})/R(\mathbf{z}^{(\varepsilon)})$ and $\delta > 0$ is the small constant appearing in Corollary 3.3. We further assume the following condition on the sequence of temperatures $\beta^{-1} = \beta(\varepsilon)^{-1} \downarrow 0$:*

$$(3.8) \qquad \beta(\varepsilon)^{-1} = o(\{\lambda^{(1,\varepsilon)} c(\varepsilon)^2 N(\varepsilon)^{-1}\}^{p/(p-1)} \lambda^{(2,\varepsilon)} \varepsilon^{\kappa/(p-1)}),$$

*as $\varepsilon \downarrow 0$. Then we have, for every $t > 0$,*

$$\lim_{\varepsilon \downarrow 0} P(\sigma^{(\varepsilon)} \ge t) = 1.$$

Theorem 3.4 asserts that asymptotically, with probability one, $\mathbf{x}^{(\varepsilon)}(t)$ keeps its rigidly crystallized shape within fluctuations $c(\varepsilon)$. In order to make the fluctuations smaller, we need better assumptions on initial data as required in (3.7) and on the speed of convergence to 0 of $\beta(\varepsilon)^{-1}$ as in (3.8). This theorem characterizes the microscopic structure of the solutions of the SDE (1.1), which are scaled macroscopically in time.

REMARK 3.2. (i) Condition $c(\varepsilon) \le \bar{c}^{(\varepsilon)}$ in (3.6) is necessary only for the set $\mathcal{M}_\infty^{\nabla,(\varepsilon)}(c(\varepsilon))$ to be well defined; recall Lemma 2.4.

(ii) Condition (3.7) is always satisfied if $\mathbf{x}^{(\varepsilon)}(0) \in \mathcal{M}^{(\varepsilon)}$.

(iii) The theorem covers the situation purely microscopic in space, that is, the case where the particles' number $N$ is fixed and does not change with $\varepsilon$. In this case, $c = c(\varepsilon)$ can be taken independently of $\varepsilon$ if it is sufficiently small and the condition (3.8) is satisfied if the temperature behaves as $\beta(\varepsilon)^{-1} = o(\varepsilon^{\kappa/(p-1)})$.

(iv) The result will be reformulated in one dimension in [6], Theorem 2.2, and the condition (3.8) will be rewritten into much simpler form based on Remarks 2.2 and 3.1 on $\lambda^{(1)}$ and $\lambda^{(2)}$, respectively.

For the proof of the theorem, we first note that $\mathbf{x}^{(\varepsilon)}(t) = (x_i^{(\varepsilon)}(t))_{i=1}^N$ satisfies the following SDE:

$$(3.9) \quad dx_i^{(\varepsilon)}(t) = -\frac{\beta}{2} \varepsilon^{-\kappa} \nabla_{x_i} H(\mathbf{x}^{(\varepsilon)}(t)) \, dt + \varepsilon^{-\kappa/2} \, dw_i(t), \qquad 1 \le i \le N,$$

in law's sense. Simple application of Itô's formula shows the next lemma.

LEMMA 3.5. *For every $p \geq 1$,*

$$H^p(\mathbf{x}^{(\varepsilon)}(t)) = H^p(\mathbf{x}^{(\varepsilon)}(0)) + m_p^{(\varepsilon)}(t)$$
$$+ \int_0^t \{-\beta\varepsilon^{-\kappa}b_{1,p}(\mathbf{x}^{(\varepsilon)}(s)) + \varepsilon^{-\kappa}b_{2,p}(\mathbf{x}^{(\varepsilon)}(s))\}\,ds,$$

*where*

$$b_{1,p}(\mathbf{x}) = \frac{p}{2}H^{p-1}(\mathbf{x})G(\mathbf{x}),$$

$$b_{2,p}(\mathbf{x}) = \frac{p}{2}(p-1)H^{p-2}(\mathbf{x})G(\mathbf{x}) + \frac{p}{2}H^{p-1}(\mathbf{x})\sum_{i=1}^N \Delta_{x_i}H(\mathbf{x}),$$

*with the Laplacian $\Delta_{x_i}$ in the variable $x_i$, and $m_p^{(\varepsilon)}(t)$ is the martingale defined by*

$$m_p^{(\varepsilon)}(t) = \varepsilon^{-\kappa/2}p\sum_{i=1}^N \int_0^t H^{p-1}(\mathbf{x}^{(\varepsilon)}(s))(\nabla_{x_i}H(\mathbf{x}^{(\varepsilon)}(s)), dw_i(s)).$$

We have the following bounds on the drift functions $b_{1,p}$ and $b_{2,p}$.

LEMMA 3.6. *Assume $\mathbf{x} \in \mathcal{M}_\infty^{\nabla,(\varepsilon)}(\delta\lambda^{(1,\varepsilon)}\lambda^{(2,\varepsilon)} \wedge \bar{c}^{(\varepsilon)})$, where $\delta > 0$ is the constant appearing in Corollary 3.3. Then, there exists $C = C_p > 0$ such that*

(3.10) $$b_{1,p}(\mathbf{x}) \geq C^{-1}\lambda^{(2,\varepsilon)}H^p(\mathbf{x}),$$

(3.11) $$b_{2,p}(\mathbf{x}) \leq CNH^{p-1}(\mathbf{x}).$$

*In particular,*

(3.12) $$-\beta b_{1,p}(\mathbf{x}) + b_{2,p}(\mathbf{x}) \leq C\beta^{-p+1}N^p\{\lambda^{(2,\varepsilon)}\}^{-p+1}.$$

PROOF. The lower bound (3.10) is immediate from Corollary 3.3. To show (3.11), note that, for every $0 < c_1 < c_2$, $|\sum_{i=1}^N \Delta_{x_i}H(\mathbf{x})| \leq CN$ holds for some $C > 0$ and all $\mathbf{x} \in (\mathbb{R}^d)^N$ satisfying $c_1 \leq |x_i - x_j| \leq c_2$ for each $\langle i, j \rangle$. Then, we have (3.11) from Corollary 3.3. Finally, to show (3.12), we estimate choosing $q = p/(p-1)$ (when $p \neq 1$),

$$H^{p-1}(\mathbf{x}) = (L\beta^{-1}N(\lambda^{(2,\varepsilon)})^{-1})^{1/q}(L^{-1}\beta N^{-1}\lambda^{(2,\varepsilon)}H^p(\mathbf{x}))^{1/q}$$
$$\leq \frac{1}{p}(L\beta^{-1}N(\lambda^{(2,\varepsilon)})^{-1})^{p-1} + \frac{1}{q}L^{-1}\beta N^{-1}\lambda^{(2,\varepsilon)}H^p(\mathbf{x}),$$

for every $L > 0$. The inequality (3.12) follows from (3.10) and (3.11) by choosing $L > 0$ sufficiently large. □

PROOF OF THEOREM 3.4. Lemmas 3.5 and 3.6 show

$$E[H^p(\mathbf{x}^{(\varepsilon)}(t \wedge \sigma))] \leq E[H^p(\mathbf{x}^{(\varepsilon)}(0))] + C\varepsilon^{-\kappa}\beta^{-p+1}N^p\{\lambda^{(2,\varepsilon)}\}^{-p+1}t.$$

However, (3.3), (2.5) and the assumption on $\mathbf{x}^{(\varepsilon)}(0)$ imply that $a_p^{(\varepsilon)} := (\lambda^{(1,\varepsilon)} \times c(\varepsilon)^2)^{-p} \times E[H^p(\mathbf{x}^{(\varepsilon)}(0))]$ tends to 0 as $\varepsilon \downarrow 0$. On the other hand, if $\mathbf{x} = (x_i)_{i=1}^N \in \partial \mathcal{M}_\infty^{\nabla,(\varepsilon)}(c(\varepsilon))$, then $|h_i - h_j| = c(\varepsilon)$ for some neighboring pair $\langle i,j \rangle$ and therefore $H(\mathbf{x}) \geq C^{-1}\lambda^{(1,\varepsilon)}c(\varepsilon)^2$ from (3.3) and (2.5). Accordingly we have

$$E[H^p(\mathbf{x}^{(\varepsilon)}(t \wedge \sigma))] \geq E[H^p(\mathbf{x}^{(\varepsilon)}(\sigma)), \sigma \leq t] \geq (C^{-1}\lambda^{(1,\varepsilon)}c(\varepsilon)^2)^p P(\sigma \leq t).$$

Therefore, we have

$$P(\sigma \leq t) \leq C\{a_p^{(\varepsilon)} + (\lambda^{(1,\varepsilon)}c(\varepsilon)^2)^{-p}\varepsilon^{-\kappa}\beta^{-p+1}N^p\{\lambda^{(2,\varepsilon)}\}^{-p+1}t\},$$

which tends to 0 as $\varepsilon \downarrow 0$. The constants $C$ may change from line to line. □

**4. Motion of a macroscopic body.** In this section we shall identify the motion of a macroscopic body obtained in the limit under the spatial scaling $x \mapsto \varepsilon x$ as $\varepsilon \downarrow 0$.

4.1. *Coordinate $\theta(\mathbf{x})$.* Let $\mathbf{z} = (z_i)_{i=1}^N$ be a centered infinitesimally rigid crystal and we fix it throughout Section 4.1. For $\mathbf{x} \in \mathcal{M}_\infty(\bar{c}(\mathbf{z}))$, $\mathbf{z}(\mathbf{x}) \in \mathcal{M}$ is defined in Section 2.2 as the minimizer of $\|\mathbf{x} - \mathbf{y}\|_2$ in $\mathbf{y} \in \mathcal{M}$; see also Section 2.5. Since $\mathbf{z}(\mathbf{x}) \in \mathcal{M}$, one can represent it as $\mathbf{z}(\mathbf{x}) = \varphi_{\theta,\eta}(\mathbf{z})$ for some $(\theta, \eta) = (\theta(\mathbf{x}), \eta(\mathbf{x})) \in SO(d) \times \mathbb{R}^d$. The function $\eta(\mathbf{x})$ defined in this way actually coincides with the center of the mass of $\mathbf{x}$:

$$(4.1) \qquad \eta(\mathbf{x}) = \frac{1}{N}\sum_{i=1}^N x_i \in \mathbb{R}^d.$$

In fact, (4.1) is seen from $d\|\mathbf{x} - \varphi_{\theta,\eta}(\mathbf{z})\|_2^2/d\eta^\alpha = 0$ for every $1 \leq \alpha \leq d$. On the other hand, the function $\theta(\mathbf{x}) = (\theta^{\alpha\beta}(\mathbf{x}))_{\alpha,\beta=1}^d$ has the following property.

LEMMA 4.1. *For every $1 \leq \alpha, \beta \leq d$ and $\mathbf{x} \in \mathcal{M}_\infty(\bar{c}(\mathbf{z}))$,*

$$(4.2) \qquad (\nabla\theta^{\alpha\beta}(\mathbf{x}), \nabla H(\mathbf{x})) \equiv \sum_{i=1}^N (\nabla_{x_i}\theta^{\alpha\beta}(\mathbf{x}), \nabla_{x_i}H(\mathbf{x})) = 0.$$

PROOF. *Step* 1. For every $\mathbf{y} \in \mathcal{M}$, let $e_\ell(\mathbf{y}) \in (\mathbb{R}^d)^N$ and $\lambda_\ell(\mathbf{y}) \geq 0$, $1 \leq \ell \leq dN$, be the eigenvectors and the corresponding eigenvalues of the Hessian $\text{Hess}\, H(\mathbf{y})$ of $H$ at $\mathbf{y}$, respectively. Recalling that $\mathbf{y}$ is infinitesimally rigid, we may assume $\lambda_\ell(\mathbf{y}) = 0$ for $1 \leq \ell \leq \tilde{d}$ and $\lambda_\ell(\mathbf{y}) > 0$ for $\tilde{d} + 1 \leq \ell \leq dN$, where $\tilde{d} := d(d+1)/2$ is the dimension of the space $\mathcal{M}$ or $\mathcal{H}_\mathbf{y}$. Then, the vectors $(e_\ell(\mathbf{y}))_{\ell=1}^{\tilde{d}}$ and $(e_\ell(\mathbf{y}))_{\ell=\tilde{d}+1}^{dN}$ span the spaces $\mathcal{H}_\mathbf{y}$ and $\mathcal{H}_\mathbf{y}^\perp$, respectively. Moreover, since the invariance of the Hamiltonian $H$: $H(\varphi_{\theta,\eta}(\mathbf{x})) = H(\mathbf{x})$ implies $\operatorname{Hess} H(\varphi_{\theta,\eta}(\mathbf{y})) = \varphi_{\theta,0} H(\mathbf{y}) \varphi_{\theta,0}^{-1}$, one can take $\{e_\ell(\mathbf{y}), \lambda_\ell(\mathbf{y})\}$ in such a manner that

$$(4.3) \qquad e_\ell(\varphi_{\theta,\eta}(\mathbf{y})) = \varphi_{\theta,0}(e_\ell(\mathbf{y})), \qquad \lambda_\ell(\varphi_{\theta,\eta}(\mathbf{y})) = \lambda_\ell(\mathbf{y}),$$

for every $1 \le \ell \le dN, (\theta, \eta) \in SO(d) \times \mathbb{R}^d$ and $\mathbf{y} \in \mathcal{M}$.

Since $\mathbf{x} - \mathbf{z}(\mathbf{x}) \in \mathcal{H}_{\mathbf{z}(\mathbf{x})}^\perp$ by (2.4), setting $\mathbf{y} := \mathbf{z}(\mathbf{x})$, $\mathbf{x} \in \mathcal{M}_\infty(\bar{c}(\mathbf{z}))$ can be decomposed as

$$(4.4) \qquad \mathbf{x} = \mathbf{y} + \sum_{\ell=\tilde{d}+1}^{dN} \zeta^\ell e_\ell(\mathbf{y})$$

for some $\zeta^\ell \in \mathbb{R}$, $\tilde{d}+1 \le \ell \le dN$. We call $(\mathbf{y}, \zeta^{\tilde{d}+1}, \ldots, \zeta^{dN}) \in \mathcal{M} \times \mathbb{R}^{dN-\tilde{d}}$ the Fermi coordinate of $\mathbf{x}$; see [7], page 4.

*Step* 2. Under the Fermi coordinate, the Hamiltonian $H$ does not depend on the variable $\mathbf{y}$:

$$(4.5) \qquad H(\mathbf{x}) = H(\zeta^{\tilde{d}+1}, \ldots, \zeta^{dN}).$$

To see (4.5), we first note that the Fermi coordinate of $\varphi_{\theta,\eta}(\mathbf{x})$ is given by $(\varphi_{\theta,\eta}(\mathbf{y}), \zeta^{\tilde{d}+1}, \ldots, \zeta^{dN})$. In fact, by (4.4) and then by (4.3),

$$\varphi_{\theta,\eta}(\mathbf{x}) = \varphi_{\theta,\eta}(\mathbf{y}) + \varphi_{\theta,0}\left(\sum_{\ell=\tilde{d}+1}^{dN} \zeta^\ell e_\ell(\mathbf{y})\right) = \varphi_{\theta,\eta}(\mathbf{y}) + \sum_{\ell=\tilde{d}+1}^{dN} \zeta^\ell e_\ell(\varphi_{\theta,\eta}(\mathbf{y})).$$

Therefore, the invariance of $H$ implies

$$H(\mathbf{y}, \zeta^{\tilde{d}+1}, \ldots, \zeta^{dN}) = H(\varphi_{\theta,\eta}(\mathbf{y}), \zeta^{\tilde{d}+1}, \ldots, \zeta^{dN})$$

under the Fermi coordinate for every $(\theta, \eta) \in SO(d) \times \mathbb{R}^d$. This shows (4.5).

*Step* 3. We finally prove $\nabla H(\mathbf{x}) \equiv (\nabla_{x_i} H(\mathbf{x}))_{i=1}^N \in \mathcal{H}_{\mathbf{z}(\mathbf{x})}^\perp$ and $\nabla \theta^{\alpha\beta}(\mathbf{x}) \in \mathcal{H}_{\mathbf{z}(\mathbf{x})}$. Once these relations are established, the conclusion of the lemma is immediately deduced. Take an arbitrary $\boldsymbol{\xi} \in (\mathbb{R}^d)^N$ and decompose it as $\boldsymbol{\xi} = P\boldsymbol{\xi} + P^\perp \boldsymbol{\xi}$, where $P : (\mathbb{R}^d)^N \to \mathcal{H}_{\mathbf{z}(\mathbf{x})}$ and $P^\perp : (\mathbb{R}^d)^N \to \mathcal{H}_{\mathbf{z}(\mathbf{x})}^\perp$ are orthogonal projections. Then

$$(4.6) \qquad \begin{aligned} (\nabla H(\mathbf{x}), \boldsymbol{\xi}) &= \frac{d}{du} H(\mathbf{x} + u\boldsymbol{\xi})\bigg|_{u=0} \\ &= \frac{d}{du} H(\mathbf{x} + uP^\perp \boldsymbol{\xi})\bigg|_{u=0} = (\nabla H(\mathbf{x}), P^\perp \boldsymbol{\xi}). \end{aligned}$$

The second equality of (4.6) follows from (4.5) noting that $\operatorname{dist}(\mathbf{z}(\mathbf{x}) + uP\boldsymbol{\xi}, \mathcal{M}) = O(u^2)$ as $u \to 0$. Equation (4.6) implies $\nabla H(\mathbf{x}) \in \mathcal{H}_{\mathbf{z}(\mathbf{x})}^\perp$. Since by definition $\theta(\mathbf{x})$ depends only on $\mathbf{y} = \mathbf{z}(\mathbf{x})$: $\theta(\mathbf{x}) = \theta(\mathbf{z}(\mathbf{x}))$, one can similarly show $\nabla \theta^{\alpha\beta}(\mathbf{x}) \in \mathcal{H}_{\mathbf{z}(\mathbf{x})}$. $\square$

In order to identify the motion of the macroscopic body, it becomes necessary to calculate the derivatives of $\theta(\mathbf{x})$ in the variables $x_i$. We introduce notation to give their representations. The space $M(d)$ of $d \times d$ matrices is equipped with an inner product $(X, Y) := \mathrm{Tr}(X^t Y)$ and a norm $|X| = \sqrt{(X, X)}$ for $X, Y \in M(d)$. The orthogonal projection from $M(d)$ onto its subspace $\mathfrak{so}(d)$ under this inner product is denoted by $\mathrm{Proj}$, that is, $\mathrm{Proj}\, X = (X - {}^t X)/2$ for $X \in M(d)$. For $e = (e^\alpha), \tilde{e} = (\tilde{e}^\alpha) \in \mathbb{R}^d$, determine the matrix $e \otimes \tilde{e} = ((e \otimes \tilde{e})^{\alpha\beta}) \in M(d)$ by $(e \otimes \tilde{e})^{\alpha\beta} = e^\alpha \tilde{e}^\beta$, $1 \leq \alpha, \beta \leq d$. The $\gamma$-directed unit vector in $\mathbb{R}^d$ is denoted by $\mathfrak{e}_\gamma, 1 \leq \gamma \leq d$. We define $Q(\mathbf{x}) = (q^{\alpha\beta}(\mathbf{x})) \in M(d)$ by

$$(4.7) \quad Q(\mathbf{x}) = \sum_{i=1}^{N} z_i \otimes x_i, \quad \text{that is,} \quad q^{\alpha\beta}(\mathbf{x}) = \sum_{i=1}^{N} z_i^\alpha x_i^\beta, \quad 1 \leq \alpha, \beta \leq d.$$

Note that $Q(\mathbf{z})$ is symmetric. The map $\{\mathrm{Proj} \circ (Q(\mathbf{x})\theta(\mathbf{x}))\}^{-1}$ is the inverse of $\mathrm{Proj} \circ (Q(\mathbf{x})\theta(\mathbf{x})) : \mathfrak{so}(d) \ni X \mapsto \mathrm{Proj}(Q(\mathbf{x})\theta(\mathbf{x})X) = \{Q(\mathbf{x})\theta(\mathbf{x})X + X^t(Q(\mathbf{x}) \times \theta(\mathbf{x}))\}/2 \in \mathfrak{so}(d)$. Note that the derivative

$$\frac{\partial \theta}{\partial x_i^\gamma}(\mathbf{x}) \in T_{\theta(\mathbf{x})}(SO(d)) \equiv \theta(\mathbf{x})\{\mathfrak{so}(d)\},$$

the tangent space to $SO(d)$ at $\theta(\mathbf{x})$.

PROPOSITION 4.2. *For $1 \leq \gamma \leq d$ and $1 \leq i \leq N$, we have*

$$\frac{\partial \theta}{\partial x_i^\gamma}(\mathbf{x}) = \theta(\mathbf{x})\{\mathrm{Proj} \circ (Q(\mathbf{x})\theta(\mathbf{x}))\}^{-1} \mathrm{Proj}\, \{(\theta(\mathbf{x})^{-1}\mathfrak{e}_\gamma) \otimes z_i\}.$$

PROOF. Since $\varphi_{\theta(\mathbf{x}), \eta(\mathbf{x})}(\mathbf{z})$ is the minimizer for the norm $\|\mathbf{x} - \mathbf{y}\|_2^2$ in $\mathbf{y} \in \mathcal{M}$, we have

$$\frac{d}{du} \sum_{j=1}^{N} |x_j - \theta(\mathbf{x}) e^{uY} z_j - \eta(\mathbf{x})|^2 \bigg|_{u=0} = 0$$

for every $Y \in \mathfrak{so}(d)$. However, since $\sum_{j=1}^{N} z_j = 0$ and $(\theta(\mathbf{x})z_j, \theta(\mathbf{x})Yz_j) = (z_j, Yz_j) = 0$ for all $j$, this implies $\sum_{j=1}^{N} (x_j, \theta(\mathbf{x})Yz_j) = 0$, which can be rewritten as

$$(4.8) \quad (Q(\mathbf{x})\theta(\mathbf{x}), Y) = 0.$$

Taking the derivative of (4.8) in $x_i^\gamma$ and noting that $\partial q^{\alpha\beta}/\partial x_i^\gamma(\mathbf{x}) = z_i^\alpha \delta^{\beta\gamma}$, we get

$$(4.9) \quad \left(Q(\mathbf{x})\frac{\partial \theta}{\partial x_i^\gamma}(\mathbf{x}), Y\right) = -\left(\frac{\partial Q}{\partial x_i^\gamma}(\mathbf{x})\theta(\mathbf{x}), Y\right) = ((\theta(\mathbf{x})^{-1}\mathfrak{e}_\gamma) \otimes z_i, Y)$$

for every $Y \in \mathfrak{so}(d)$. Set $X := \theta(\mathbf{x})^{-1} \partial \theta / \partial x_i^\gamma(\mathbf{x})$ and then, since $X \in \mathfrak{so}(d)$, (4.9) shows that

$$\text{Proj}\,(Q(\mathbf{x})\theta(\mathbf{x})X) = \text{Proj}\,\{(\theta(\mathbf{x})^{-1}\mathfrak{e}_\gamma) \otimes z_i\}.$$

This proves the conclusion. $\square$

4.2. *Identification of the limit.* We now discuss the limit as $\varepsilon \downarrow 0$ under the macroscopic spatial scaling $x \mapsto y = \varepsilon x$ for the system with the particles' number $N \equiv N(\varepsilon)$ changing with $\varepsilon$.

Our formulation is the following. Let $\mathfrak{M}_{\tilde{\rho}} \equiv \mathfrak{M}_{\tilde{\rho}}(\mathbb{R}^d)$, $\tilde{\rho} > 0$, be the family of all Radon measures $\mu$ on $\mathbb{R}^d$ satisfying $\mu(\mathbb{R}^d) \leq \tilde{\rho}$. The space $\mathfrak{M}_{\tilde{\rho}}$ is equipped with the topology determined by the weak convergence. A sequence $\mathbf{x}^{(\varepsilon)} = (x_i^{(\varepsilon)})_{i=1}^N$, $0 < \varepsilon < 1$, with $N = N(\varepsilon)$ of the system of particles in $\mathbb{R}^d$ is identified under the scaling with $\mu^{(\varepsilon)}(\mathbf{x}^{(\varepsilon)}) \in \mathfrak{M}_{\tilde{\rho}}, \tilde{\rho} = \varepsilon^d N$ defined by

$$(4.10) \qquad \mu^{(\varepsilon)}(dy) \equiv \mu^{(\varepsilon)}(\mathbf{x}^{(\varepsilon)}; dy) := \varepsilon^d \sum_{i=1}^N \delta_{\varepsilon x_i^{(\varepsilon)}}(dy).$$

Let us assume that, as in Section 3.2, a sequence $\mathbf{z}^{(\varepsilon)} \equiv (z_i^{(\varepsilon)})_{i=1}^N$, $0 < \varepsilon < 1$, of centered infinitesimally rigid crystals is given and satisfies the following three conditions:

1. There exists $R > 0$ such that $|z_i^{(\varepsilon)}| \leq R\varepsilon^{-1}$ for all $i$ and $\varepsilon$.
2. $\mathbf{z}^{(\varepsilon)}$ has a *macroscopic limit density function* $\rho(y)$, $y \in \mathbb{R}^d$, in the sense that $\mu^{(\varepsilon)}(\mathbf{z}^{(\varepsilon)}; dy) \Longrightarrow \rho(y)\,dy$ weakly as $\varepsilon \downarrow 0$.
3. (Nontriviality of the limit). The total mass of the macroscopic limit density is positive: $\bar{\rho} := \int_{\mathbb{R}^d} \rho(y)\,dy > 0$.

Examples of the sequence $\mathbf{z}^{(\varepsilon)}$ will be given at the end of this section. We shall denote by $\mathcal{M}_\infty^{(\varepsilon)}(c)$ the set (2.7) determined from $\mathbf{z}^{(\varepsilon)}$ instead of $\mathbf{z}$. The domain $D := \{y \in \mathbb{R}^d; \rho(y) > 0\}$ with density $\rho(y)$ may be called the macroscopic shape of the body. The above conditions imply that $D \subset \{y \in \mathbb{R}^d; |y| \leq R\}$ and $\lim_{\varepsilon \downarrow 0} \varepsilon^d N(\varepsilon) = \bar{\rho}$. Since $\mathbf{z}^{(\varepsilon)}$ are centered, the body with density $\rho(y)$ is also centered in the sense that $\int_{\mathbb{R}^d} y\rho(y)\,dy = 0$. Let $\bar{Q} = (\bar{q}^{\alpha\beta})_{\alpha,\beta=1}^d \in M(d)$ be the matrix defined by

$$\bar{q}^{\alpha\beta} = \int_{\mathbb{R}^d} y^\alpha y^\beta \rho(y)\,dy.$$

Then, since the matrix $\bar{Q}$ is symmetric, by rotating the body around the origin $0 \in \mathbb{R}^d$, we may assume that $\bar{Q}$ is diagonal with diagonal elements:

$$(4.11) \qquad \bar{q}^\alpha = \int_{\mathbb{R}^d} (y^\alpha)^2 \rho(y)\,dy, \qquad 1 \leq \alpha \leq d.$$

Let $\mathbf{x}^{(\varepsilon)}(t) := \mathbf{x}(\varepsilon^{-\kappa}t)$, $\kappa = d+2$, be the process obtained by macroscopically scaling in time the solution $\mathbf{x}(t)$ of the SDE (1.1) with initial configuration $\mathbf{x}(0) = \mathbf{z}^{(\varepsilon)}$. The spatially macroscopic scaling limit of $\mathbf{x}^{(\varepsilon)}(t)$ is characterized by the following theorem, in which the limit of $(\eta^{(\varepsilon)}(t), \theta^{(\varepsilon)}(t)) := (\varepsilon\eta(\mathbf{x}^{(\varepsilon)}(t)), \theta(\mathbf{x}^{(\varepsilon)}(t))) \in \mathbb{R}^d \times SO(d)$ as $\varepsilon \downarrow 0$ is obtained. Here, the coordinate $\theta(\mathbf{x})$, $\mathbf{x} \in \mathcal{M}_\infty^{(\varepsilon)}(\bar{c}(\mathbf{z}^{(\varepsilon)}))$, is defined as in Section 4.1 based on $\mathbf{z}^{(\varepsilon)}$ in place of $\mathbf{z}$. If $\mathbf{x}^{(\varepsilon)}(t)$ goes outside of $\mathcal{M}_\infty^{(\varepsilon)}(\bar{c}(\mathbf{z}^{(\varepsilon)}))$ at a certain time, $\theta^{(\varepsilon)}(t)$ may be defined arbitrarily after such time keeping it continuous in $t$. This theorem, in particular, shows that the proper macroscopic time scalings for the translational and rotational motions of the body are the same.

THEOREM 4.3. *Assume that the temperature $\beta^{-1} = \beta(\varepsilon)^{-1}$ of the system converges to $0$ as in (3.8) for some $c(\varepsilon)$ satisfying (3.6) and $c(\varepsilon) \leq \varepsilon^\nu$ for some $\nu > 0$. Then, the process $(\eta^{(\varepsilon)}(t), \theta^{(\varepsilon)}(t))$ weakly converges to $(\eta(t), \theta(t))$ as $\varepsilon \downarrow 0$ in the space $C([0,T], \mathbb{R}^d \times SO(d))$ for every $T > 0$. The limit is characterized by the following three properties:*

 (i) *$\eta(t)$ and $\theta(t)$ are mutually independent.*
 (ii) *$\sqrt{\bar{\rho}}\eta(t)$ is a $d$-dimensional Brownian motion starting at $0$.*
 (iii) *$\theta(t)$ is a solution of an SDE of Stratonovich's type on $SO(d)$:*

(4.12) $$d\theta(t) = \theta(t) \circ dm(t), \qquad \theta(0) = I,$$

*where $m(t) = (m^{\alpha\beta}(t))_{\alpha,\beta=1}^d$ is an $\mathfrak{so}(d)$-valued Brownian motion such that the components $\{m^{\alpha\beta}(t); \alpha < \beta\}$ in the upper half of the matrix $m(t)$ are mutually independent and $\sqrt{\bar{q}^\alpha + \bar{q}^\beta}m^{\alpha\beta}(t)$ is one-dimensional Brownian motion for each $1 \leq \alpha, \beta \leq d$.*

This theorem can be reformulated as the convergence for measure-valued processes.

COROLLARY 4.4. *Under the same assumption on $\beta(\varepsilon)^{-1}$ as Theorem 4.3, $\mu^{(\varepsilon)}(t) := \mu^{(\varepsilon)}(\mathbf{x}^{(\varepsilon)}(t);\cdot)$ weakly converges to $\mu(t) := \rho(\varphi^{-1}_{\theta(t),\eta(t)}(y))\,dy$ as $\varepsilon \downarrow 0$ in the space $C([0,T], \mathfrak{M}_{\tilde\rho})$ for every $T > 0$, where $\tilde\rho = \sup_{0<\varepsilon<1}\varepsilon^d N(\varepsilon)$. The process $(\theta(t), \eta(t))$ is characterized by the three properties (i)–(iii) in Theorem 4.3 and $\rho(y)$ is the macroscopic limit density function of the initial configuration $\mathbf{z}^{(\varepsilon)}$.*

Before giving the proofs of Theorem 4.3 and Corollary 4.4, we state a proposition whose proof will be postponed to the next section.

PROPOSITION 4.5. *For every $\nu > 0$, there exist $C > 0$ and $0 < \varepsilon_0 < 1$ such that*

$$(4.13) \qquad \left| \frac{\partial \theta}{\partial x_i^\gamma}(\mathbf{x}) \right| \leq C \varepsilon^{\kappa - 1},$$

$$(4.14) \qquad \left| \frac{\partial \theta}{\partial x_i^\gamma}(\mathbf{x}) - \frac{\partial \theta}{\partial x_i^\gamma}(\mathbf{z}(\mathbf{x})) \right| \leq C \varepsilon^{d + \nu + 1},$$

*hold for every $\mathbf{x} \in \mathcal{M}_\infty^{(\varepsilon)}(\varepsilon^{\nu - 1} \wedge \bar{c}(\mathbf{z}^{(\varepsilon)}))$, $1 \leq \gamma \leq d$, $1 \leq i \leq N(\varepsilon)$ and $0 < \varepsilon < \varepsilon_0$.*

PROOF OF THEOREM 4.3. *Step* 1. From the SDE (3.9) for $\mathbf{x}^{(\varepsilon)}(t)$, since $\nabla U(x) = -\nabla U(-x)$ holds from the radial symmetry of $U$, we have

$$(4.15) \qquad \eta^{(\varepsilon)}(t) = \frac{\varepsilon^{1 - \kappa/2}}{N} \sum_{i=1}^N w_i(t),$$

which is equivalent to $\varepsilon^{1-\kappa/2} N^{-1/2} w(t)$ in law. This shows the property (ii) for the limit $\eta(t)$ of $\eta^\varepsilon(t)$ noting that $\bar{\rho} = \lim_{\varepsilon \downarrow 0} \varepsilon^d N$ so that $\lim_{\varepsilon \downarrow 0} \varepsilon^{1-\kappa/2} N^{-1/2} = \bar{\rho}^{-1/2}$.

*Step* 2. We next consider the limit of $\theta^{(\varepsilon)}(t)$. Let $\sigma^{(\varepsilon)}$ be the stopping time defined by (3.5) with $c(\varepsilon)$ satisfying the conditions in the theorem. Then, $\theta(\mathbf{x}^{(\varepsilon)}(t))$ is well defined for $t \leq \sigma^{(\varepsilon)}$, and again from the SDE (3.9) and using the property (4.2) of $\theta(\mathbf{x})$, we have

$$(4.16) \qquad d\theta^{(\varepsilon)}(t) = \varepsilon^{-\kappa/2} \sum_{i=1}^N \nabla_{x_i} \theta(\mathbf{x}^{(\varepsilon)}(t)) \circ dw_i(t), \qquad t \leq \sigma^{(\varepsilon)},$$

by applying Itô's formula, where $\nabla_{x_i} \theta \circ dw_i := \sum_{\gamma=1}^d \partial \theta / \partial x_i^\gamma \circ dw_i^\gamma \in M(d)$. This may be further rewritten as

$$d\theta^{(\varepsilon)}(t) = \theta^{(\varepsilon)}(t) \circ dm^{(\varepsilon)}(t), \qquad t \leq \sigma^{(\varepsilon)},$$

with an $\mathfrak{so}(d)$-valued martingale $m^{(\varepsilon)}(t) = \{m^{\alpha\beta,(\varepsilon)}(t)\}_{\alpha,\beta=1}^d$ defined by

$$m^{(\varepsilon)}(t) = \varepsilon^{-\kappa/2} \sum_{i=1}^N \int_0^t \theta(\mathbf{x}^{(\varepsilon)}(s))^{-1} \nabla_{x_i} \theta(\mathbf{x}^{(\varepsilon)}(s)) \, dw_i(s), \qquad t \leq \sigma^{(\varepsilon)}.$$

Therefore, setting $\tilde{\theta}^{(\varepsilon)}(t) := \theta^{(\varepsilon)}(t \wedge \sigma^{(\varepsilon)})$ and $\tilde{m}^{(\varepsilon)}(t) := m^{(\varepsilon)}(t \wedge \sigma^{(\varepsilon)})$ for all $t \geq 0$, $\tilde{\theta}^{(\varepsilon)}(t)$ satisfies an SDE

$$d\tilde{\theta}^{(\varepsilon)}(t) = \tilde{\theta}^{(\varepsilon)}(t) \circ d\tilde{m}^{(\varepsilon)}(t), \qquad t \geq 0.$$

This SDE written in Stratonovich's form is equivalent to

$$d\tilde{\theta}^{(\varepsilon)}(t) = \tilde{\theta}^{(\varepsilon)}(t) d\tilde{m}^{(\varepsilon)}(t) + \tfrac{1}{2} \tilde{\theta}^{(\varepsilon)}(t) d\langle \tilde{m}^{(\varepsilon)}, \tilde{m}^{(\varepsilon)} \rangle(t)$$

in Itô's form, where the quadratic variational process $\langle \tilde{m}^{(\varepsilon)}, \tilde{m}^{(\varepsilon)} \rangle(t) \in M(d)$ is defined by

$$(\langle \tilde{m}^{(\varepsilon)}, \tilde{m}^{(\varepsilon)} \rangle(t))^{\alpha\beta} = \sum_{\gamma=1}^{d} \langle \tilde{m}^{\alpha\gamma,(\varepsilon)}, \tilde{m}^{\gamma\beta,(\varepsilon)} \rangle(t), \qquad 1 \leq \alpha, \beta \leq d.$$

The goal is to show that $\tilde{\theta}^{(\varepsilon)}(t)$ weakly converges as $\varepsilon \downarrow 0$ to the solution $\theta(t)$ of the SDE (4.12) which is equivalent to

$$d\theta(t) = \theta(t)dm(t) + \tfrac{1}{2}\theta(t)d\langle m, m \rangle(t)$$

in Itô's form. Indeed, once this is proved, since $\lim_{\varepsilon \downarrow 0} P(\sigma^{(\varepsilon)} \geq T) = 1$ from the microscopic shape theorem (Theorem 3.4), $\theta^{(\varepsilon)}(t)$ also weakly converges to $\theta(t)$; recall that $\theta^{(\varepsilon)}(t)$ was arbitrarily defined after the time when $\mathbf{x}^{(\varepsilon)}(t)$ goes outside of $\mathcal{M}_{\infty}^{(\varepsilon)}(\bar{c}(\mathbf{z}^{(\varepsilon)}))$. To show the weak convergence of $\tilde{\theta}^{(\varepsilon)}(t)$ to $\theta(t)$ in $C([0,T], SO(d))$, it suffices to prove the following two conditions for the driving martingale $\tilde{m}^{(\varepsilon)}(t)$:

$$(4.17) \quad \lim_{\varepsilon \downarrow 0} E\left[\left| E[\langle \tilde{m}^{\alpha\beta,(\varepsilon)}, \tilde{m}^{ab,(\varepsilon)} \rangle(t) | \mathcal{F}_s] - \frac{t(\delta^{\alpha a}\delta^{\beta b} - \delta^{\alpha b}\delta^{\beta a})}{\bar{q}^{\alpha} + \bar{q}^{\beta}} \right|\right] = 0$$

for every $t \geq s \geq 0$, $1 \leq \alpha, \beta, a, b \leq d$, where $\mathcal{F}_s = \sigma\{w(s'); s' \leq s\}$ and

$$(4.18) \qquad \sup_{0<\varepsilon<1} \sup_{0 \leq t \leq T} E\left[\left|\frac{d}{dt}\langle \tilde{m}^{\alpha\beta,(\varepsilon)}, \tilde{m}^{ab,(\varepsilon)} \rangle(t)\right|^p\right] < \infty$$

for some $p > 2$; see, for instance, [11], page 222, Theorem 5.2.1. Note that, from the property (iii), the quadratic variational processes of $(m^{\alpha\beta}(t))_{\alpha\beta}$ are given by

$$\langle m^{\alpha\beta}, m^{ab} \rangle(t) = \frac{t(\delta^{\alpha a}\delta^{\beta b} - \delta^{\alpha b}\delta^{\beta a})}{\bar{q}^{\alpha} + \bar{q}^{\beta}}.$$

Now let us prove (4.17) and (4.18). Since we have

$$\frac{d}{dt}\langle \tilde{m}^{\alpha\beta,(\varepsilon)}, \tilde{m}^{ab,(\varepsilon)} \rangle(t)$$
$$= \varepsilon^{-\kappa} \sum_{i=1}^{N} \sum_{\gamma=1}^{d} \left(\theta(\mathbf{x}^{(\varepsilon)}(t))^{-1} \frac{\partial \theta}{\partial x_i^{\gamma}}(\mathbf{x}^{(\varepsilon)}(t))\right)^{\alpha\beta} \left(\theta(\mathbf{x}^{(\varepsilon)}(t))^{-1} \frac{\partial \theta}{\partial x_i^{\gamma}}(\mathbf{x}^{(\varepsilon)}(t))\right)^{ab}$$

for $t \leq \sigma^{(\varepsilon)}$ (and $= 0$ for $t > \sigma^{(\varepsilon)}$), (4.13) in Proposition 4.5 implies

$$(4.19) \qquad \left|\frac{d}{dt}\langle \tilde{m}^{\alpha\beta,(\varepsilon)}, \tilde{m}^{ab,(\varepsilon)} \rangle(t)\right| \leq C\varepsilon^{-\kappa} N(\varepsilon^{\kappa-1})^2 \leq C',$$

which is bounded in $\varepsilon$ and $t$. This shows (4.18). To prove (4.17), define $\bar{m}^{(\varepsilon)}(t)$ by

$$(4.20) \quad \bar{m}^{(\varepsilon)}(t) = \varepsilon^{-\kappa/2} \sum_{i=1}^{N} \int_{0}^{t \wedge \sigma^{(\varepsilon)}} \theta(\mathbf{z}(\mathbf{x}^{(\varepsilon)}(s)))^{-1} \nabla_{x_i} \theta(\mathbf{z}(\mathbf{x}^{(\varepsilon)}(s))) \, dw_i(s).$$

We have replaced $\mathbf{x}^{(\varepsilon)}(s)$ with $\mathbf{z}(\mathbf{x}^{(\varepsilon)}(s))$ in $\tilde{m}^{(\varepsilon)}(t)$. Then, since condition 1 on $\mathbf{z}^{(\varepsilon)}$ implies $R(\mathbf{z}^{(\varepsilon)}) \le R_0 \varepsilon^{-1}$ for some $R_0 > 0$, Lemma 2.4 shows $\mathbf{x}^{(\varepsilon)}(s) \in \mathcal{M}_{\infty}^{\nabla,(\varepsilon)}(c(\varepsilon)) \subset \mathcal{M}_{\infty}^{(\varepsilon)}(R(\mathbf{z}^{(\varepsilon)})c(\varepsilon)) \subset \mathcal{M}_{\infty}^{(\varepsilon)}(R_0 \varepsilon^{\nu-1} \wedge \bar{c}(\mathbf{z}^{(\varepsilon)})) \subset \mathcal{M}_{\infty}^{(\varepsilon)}(\varepsilon^{\nu'-1} \wedge \bar{c}(\mathbf{z}^{(\varepsilon)}))$ for $s \le \sigma^{(\varepsilon)}$ by taking $\nu' \in (0, \nu)$ and for small $\varepsilon > 0$. Therefore, by (4.14) in Proposition 4.5 and noting that $\theta(\mathbf{x}) = \theta(\mathbf{z}(\mathbf{x}))$, we have

$$(4.21) \quad |\langle \tilde{m}^{\alpha\beta,(\varepsilon)} - \bar{m}^{\alpha\beta,(\varepsilon)} \rangle(t)| \le C \varepsilon^{-\kappa} N (\varepsilon^{d+\nu'+1})^2 t \le C' \varepsilon^{2\nu'} t.$$

Using (4.19) and (4.21), we get

$$|\langle \tilde{m}^{\alpha\beta,(\varepsilon)}, \tilde{m}^{ab,(\varepsilon)} \rangle(t) - \langle \bar{m}^{\alpha\beta,(\varepsilon)}, \bar{m}^{ab,(\varepsilon)} \rangle(t)|$$
$$\le \sqrt{\langle \tilde{m}^{\alpha\beta,(\varepsilon)} - \bar{m}^{\alpha\beta,(\varepsilon)} \rangle(t) \cdot \langle \tilde{m}^{ab,(\varepsilon)} \rangle(t)}$$
$$+ \sqrt{\langle \tilde{m}^{ab,(\varepsilon)} - \bar{m}^{ab,(\varepsilon)} \rangle(t) \cdot \langle \bar{m}^{\alpha\beta,(\varepsilon)} \rangle(t)}$$
$$\le C'' \varepsilon^{\nu'} t,$$

which tends to 0 as $\varepsilon \downarrow 0$; note that a similar estimate to (4.19) can be shown for $\langle \bar{m}^{\alpha\beta,(\varepsilon)}, \bar{m}^{ab,(\varepsilon)} \rangle(t)$. We accordingly see that it suffices to show condition (4.17) for $\bar{m}^{(\varepsilon)}$ instead of $\tilde{m}^{(\varepsilon)}$.

For a configuration $\tilde{\mathbf{z}} = \varphi_{\tilde{\theta},\tilde{\eta}}(\mathbf{z}^{(\varepsilon)}) \in \mathcal{M}^{(\varepsilon)}$ with some $(\tilde{\theta}, \tilde{\eta}) \in SO(d) \times \mathbb{R}^d$, since $\theta(\tilde{\mathbf{z}}) = \tilde{\theta}$ and $Q(\tilde{\mathbf{z}}) = Q(\mathbf{z}^{(\varepsilon)}) \tilde{\theta}^{-1}$, Proposition 4.2 shows that

$$\frac{\partial \theta}{\partial x_i^{\gamma}} (\varphi_{\tilde{\theta},\tilde{\eta}}(\mathbf{z}^{(\varepsilon)})) = \tilde{\theta} \{\operatorname{Proj} \circ (Q(\mathbf{z}^{(\varepsilon)}))\}^{-1} \operatorname{Proj} \{(\tilde{\theta}^{-1} \mathfrak{e}_{\gamma}) \otimes z_i^{(\varepsilon)}\}.$$

However, we see that

$$(4.22) \quad \lim_{\varepsilon \downarrow 0} \varepsilon^{\kappa} Q(\mathbf{z}^{(\varepsilon)}) = \lim_{\varepsilon \downarrow 0} \left( \int_{\mathbb{R}^d} y^{\alpha} y^{\beta} \mu^{(\varepsilon)}(\mathbf{z}^{(\varepsilon)}; dy) \right)_{\alpha,\beta=1}^{d} = \bar{Q},$$

$$(4.23) \quad \{\operatorname{Proj} \circ (\bar{Q})\}^{-1} Y = \left( \frac{2 y^{\alpha\beta}}{\bar{q}^{\alpha} + \bar{q}^{\beta}} \right)_{\alpha,\beta=1}^{d},$$

for $Y = (y^{\alpha\beta})_{\alpha\beta} \in \mathfrak{so}(d)$, and

$$\operatorname{Proj} \{(\tilde{\theta}^{-1} \mathfrak{e}_{\gamma}) \otimes z_i^{(\varepsilon)}\} = \tfrac{1}{2} (\tilde{\theta}^{\gamma\alpha} z_i^{\beta,(\varepsilon)} - \tilde{\theta}^{\gamma\beta} z_i^{\alpha,(\varepsilon)})_{\alpha,\beta=1}^{d},$$

where $\tilde{\theta} = (\tilde{\theta}^{\alpha\beta})_{\alpha\beta}$ and $z_i^{(\varepsilon)} = (z_i^{\alpha,(\varepsilon)})_{\alpha}$. Therefore, from Lemma 5.1(ii) stated below [note that $\{\operatorname{Proj} \circ (Q(\mathbf{z}^{(\varepsilon)}))\}^{-1} = \varepsilon^{\kappa} \Phi^{(\varepsilon)}(\mathbf{z}^{(\varepsilon)})^{-1}$] and recalling $\varepsilon |z_i^{(\varepsilon)}| \le$

$R$, we have

$$(4.24) \quad \lim_{\varepsilon \downarrow 0} \left| \varepsilon^{-\kappa+1} \tilde{\theta}^{-1} \frac{\partial \theta}{\partial x_i^\gamma}(\varphi_{\tilde{\theta},\tilde{\eta}}(\mathbf{z}^{(\varepsilon)})) - \varepsilon \left( \frac{\tilde{\theta}^{\gamma\alpha} z_i^{\beta,(\varepsilon)} - \tilde{\theta}^{\gamma\beta} z_i^{\alpha,(\varepsilon)}}{\bar{q}^\alpha + \bar{q}^\beta} \right)_{\alpha,\beta=1}^d \right| = 0$$

uniformly in $(\tilde{\theta}, \tilde{\eta})$. Hence,

$$\frac{d}{dt} \langle \bar{m}^{\alpha\beta,(\varepsilon)}, \bar{m}^{ab,(\varepsilon)} \rangle(t)$$

$$= \frac{\varepsilon^\kappa}{(\bar{q}^\alpha + \bar{q}^\beta)(\bar{q}^a + \bar{q}^b)} \sum_{i=1}^N \sum_{\gamma=1}^d (\theta^{\gamma\alpha,(\varepsilon)}(t) z_i^{\beta,(\varepsilon)} - \theta^{\gamma\beta,(\varepsilon)}(t) z_i^{\alpha,(\varepsilon)})$$

$$\times (\theta^{\gamma a,(\varepsilon)}(t) z_i^{b,(\varepsilon)} - \theta^{\gamma b,(\varepsilon)}(t) z_i^{a,(\varepsilon)}) + o(1)$$

$$= \frac{\delta^{\alpha a} \delta^{\beta b} - \delta^{\alpha b} \delta^{\beta a}}{\bar{q}^\alpha + \bar{q}^\beta} + o(1), \qquad t \leq \sigma^{(\varepsilon)},$$

as $\varepsilon \downarrow 0$ for every $1 \leq \alpha, \beta, a, b \leq d$. The error terms $o(1)$, which come from the errors in (4.24), tend to 0 uniformly in $t$ and $\omega$ [an element of the probability space on which $\bar{m}^{(\varepsilon)}(t)$ are defined]. This proves (4.17) for $\bar{m}^{(\varepsilon)}(t)$.

*Step* 3. Finally, to show property (i), compute the quadratic variational processes of $\bar{m}^{\alpha\beta,(\varepsilon)}(t)$ and $\eta^{\gamma,(\varepsilon)}(t)$ from (4.15), (4.20) and (4.24):

$$\frac{d}{dt} \langle \bar{m}^{\alpha\beta,(\varepsilon)}, \eta^{\gamma,(\varepsilon)} \rangle(t)$$

$$= \frac{\varepsilon}{N} \sum_{i=1}^N \frac{\theta^{\gamma\alpha,(\varepsilon)}(t) z_i^{\beta,(\varepsilon)} - \theta^{\gamma\beta,(\varepsilon)}(t) z_i^{\alpha,(\varepsilon)}}{\bar{q}^\alpha + \bar{q}^\beta} + o(1), \qquad t \leq \sigma^{(\varepsilon)}.$$

However, since $\mathbf{z}^{(\varepsilon)}$ is centered, the sum vanishes and this proves property (i). □

PROOF OF COROLLARY 4.4. Theorem 4.3 combined with Theorem 3.4 shows that $\langle f, \mu^{(\varepsilon)}(t) \rangle$ weakly converges to $\langle f, \mu(t) \rangle$ as $\varepsilon \downarrow 0$ in the space $C([0,T], \mathbb{R})$ for every $f \in C_b(\mathbb{R}^d)$, where $\langle f, \mu \rangle := \int_{\mathbb{R}^d} f(y) \mu(dy)$. Therefore, to conclude the corollary, it suffices to show the tightness of the family of laws of $\{\mu^{(\varepsilon)}(t); 0 < \varepsilon < 1\}$ on the space $C([0,T], \mathfrak{M}_{\tilde{\rho}})$. But this can be deduced from:

1. For each $\delta > 0$, there exists a compact set $\mathcal{K}$ in $\mathfrak{M}_{\tilde{\rho}}$ such that $P(\mu^{(\varepsilon)}(t) \in \mathcal{K}$ for every $t \in [0,T]) \geq 1 - \delta$.
2. For every $f \in C_b(\mathbb{R}^d)$, $\{\langle f, \mu^{(\varepsilon)}(t) \rangle; 0 < \varepsilon < 1\}$ is tight on the space $C([0,T], \mathbb{R})$.

Condition 2 follows from what we mentioned above. Condition 1 is also easy, since the support of $\mu^{(\varepsilon)}(t)$ is in the ball with center $\eta^{(\varepsilon)}(t)$ and radius $R$ as long as $t \leq \sigma^{(\varepsilon)}$ and $\eta^{(\varepsilon)}(t) \Longrightarrow \eta(t)$, which is the time changed Brownian motion, as $\varepsilon \downarrow 0$. □

The solution $\theta(t)$ of the SDE (4.12) is called the left Brownian motion on $SO(d)$; see [15, 16]. A coordinate satisfying the relation like (4.2) was used by Katzenberger [9] to make a cancellation for diverging terms as we have seen in deriving (4.16).

EXAMPLE 4.1. (i) Let $D$ be a bounded domain in $\mathbb{R}^d (d=2$ or $3)$ having a smooth boundary $\partial D$ and let $\mathbf{z}^{(\varepsilon)}$ be the infinitesimally rigid crystal constructed from $\varepsilon^{-1}D \cap a\Lambda$ on the $d$-dimensional triangular lattice as in Example 2.3. The configuration $\mathbf{z}^{(\varepsilon)}$ is the microscopic crystal consisting of atoms arranged in an equal distance $a$ and $D$ is the corresponding macroscopic body. The macroscopic density function of $\mathbf{z}^{(\varepsilon)}$ is given by

$$\rho(y) = \frac{\mathbb{1}_D(y)}{a^d |\det A|},$$

where $A = (e_1 e_2 \cdots e_d) \in M(d)$ is the matrix consisting of $d$ column vectors $\{e_\alpha\}_\alpha$ used for the definition of the triangular lattice. In this sense, $D$ is the high density region and the outside of $D$ is the empty region.

(ii) In higher dimensions, one can construct $\mathbf{z}^{(\varepsilon)}$ based on the idea explained in Example 2.1.

**5. Proof of Proposition 4.5.** This section gives the proof of Proposition 4.5. Consider an operator $\Phi$ on the space $\mathfrak{so}(d)$ defined by $\Phi X = \mathrm{Proj}(\bar{Q}X)$ for $X \in \mathfrak{so}(d)$. Then, as we have seen in (4.23), $\Phi$ is invertible and the operator norm of $\Phi^{-1}$ can be dominated by

$$(5.1) \qquad \|\Phi^{-1}\| \leq \bar{C} := \max_{1 \leq \alpha \neq \beta \leq d} \frac{2}{\bar{q}^\alpha + \bar{q}^\beta}.$$

Recall that $Q(\mathbf{x}) \equiv Q^{(\varepsilon)}(\mathbf{x})$ is determined from $\mathbf{z}^{(\varepsilon)}$ and set $\bar{Q}^{(\varepsilon)}(\mathbf{x}) = \varepsilon^\kappa Q^{(\varepsilon)}(\mathbf{x})$. Then, we have $\bar{Q} = \lim_{\varepsilon \downarrow 0} \bar{Q}^{(\varepsilon)}(\mathbf{z}^{(\varepsilon)})$; see (4.22). We introduce another operator $\Phi^{(\varepsilon)}(\mathbf{x})$ on $\mathfrak{so}(d)$ by

$$\Phi^{(\varepsilon)}(\mathbf{x})X = \mathrm{Proj}\,(\bar{Q}^{(\varepsilon)}(\mathbf{x})\theta(\mathbf{x})X), \qquad X \in \mathfrak{so}(d).$$

In the following, we denote $\mathcal{M}_\infty^{(\varepsilon)}(\varepsilon^{\nu-1} \wedge \bar{c}(\mathbf{z}^{(\varepsilon)}))$ simply by $\mathcal{M}_\infty^{(\varepsilon)}(\varepsilon^{\nu-1})$, since we only use the bound $\|\mathbf{h}(\mathbf{x})\|_\infty \leq \varepsilon^{\nu-1}$ for $\mathbf{x}$.

LEMMA 5.1. (i) *For every $\nu > 0$, there exists $\varepsilon_0 > 0$ such that*

$$\sup_{0 < \varepsilon \leq \varepsilon_0} \sup_{\mathbf{x} \in \mathcal{M}_\infty^{(\varepsilon)}(\varepsilon^{\nu-1})} \|\Phi^{(\varepsilon)}(\mathbf{x})^{-1}\| < \infty.$$

(ii) *For every $\nu > 0$,*

$$\lim_{\varepsilon \downarrow 0} \sup_{\mathbf{x} \in \mathcal{M}_\infty^{(\varepsilon)}(\varepsilon^{\nu-1})} \|\Phi^{(\varepsilon)}(\mathbf{x})^{-1} - \Phi^{-1}\| = 0.$$

PROOF. We first assume that $\mathbf{x} \in \mathcal{M}_\infty^{(\varepsilon)}(\varepsilon^{\nu-1})$ satisfies $\mathbf{z}(\mathbf{x}) = \mathbf{z}^{(\varepsilon)}$. Then, we have

(5.2) $$|\bar{Q}^{(\varepsilon)}(\mathbf{x}) - \bar{Q}| \leq r(\varepsilon),$$

with some $r(\varepsilon) \to 0$ as $\varepsilon \downarrow 0$. Indeed, the left-hand side of (5.2) is dominated by

$$|\bar{Q}^{(\varepsilon)}(\mathbf{x}) - \bar{Q}^{(\varepsilon)}(\mathbf{z}^{(\varepsilon)})| + |\bar{Q}^{(\varepsilon)}(\mathbf{z}^{(\varepsilon)}) - \bar{Q}|,$$

and the first term is further bounded as

$$|\varepsilon^\kappa q^{\alpha\beta}(\mathbf{x}) - \varepsilon^\kappa q^{\alpha\beta}(\mathbf{z}^{(\varepsilon)})| \leq \varepsilon^\kappa \sum_{i=1}^N |z_i^{\alpha,(\varepsilon)}||x_i^\beta - z_i^{\beta,(\varepsilon)}|$$

$$\leq \varepsilon^\kappa \cdot N \cdot R\varepsilon^{-1}\varepsilon^{\nu-1} \leq C\varepsilon^\nu,$$

while the second term tends to 0 as $\varepsilon \downarrow 0$.

Denoting by $\Psi = \Phi^{(\varepsilon)}(\mathbf{x}) - \Phi$, $\Phi^{(\varepsilon)}(\mathbf{x})^{-1}$ can be expressed as

$$\Phi^{(\varepsilon)}(\mathbf{x})^{-1} = (I + \Phi^{-1}\Psi)^{-1}\Phi^{-1} = \sum_{k=0}^\infty (-\Phi^{-1}\Psi)^k \Phi^{-1}.$$

Since $\theta(\mathbf{x}) = \theta(\mathbf{z}^{(\varepsilon)}) = I$, we have $\Psi X = \frac{1}{2}\{\bar{Q}^{(\varepsilon)}(\mathbf{x}) - \bar{Q}\}X + \frac{1}{2}X^t\{\bar{Q}^{(\varepsilon)}(\mathbf{x}) - \bar{Q}\}$ and therefore $\|\Psi\| \leq r(\varepsilon)$ from (5.2). Accordingly, (5.1) implies

$$\left\|\sum_{k=0}^\infty (-\Phi^{-1}\Psi)^k\right\| \leq \frac{1}{1 - \bar{C}r(\varepsilon)} \leq 2,$$

so $\|\Phi^{(\varepsilon)}(\mathbf{x})^{-1}\| \leq 2\bar{C}$ for sufficiently small $\varepsilon > 0$ such that $r(\varepsilon) \leq 1/(2\bar{C})$. By acting rotation and translation, similar estimate can be derived for $\|\Phi^{(\varepsilon)}(\mathbf{x})^{-1}\|$ for every $\mathbf{x} \in \mathcal{M}_\infty^{(\varepsilon)}(\varepsilon^{\nu-1})$ [without assuming $\mathbf{z}(\mathbf{x}) = \mathbf{z}^{(\varepsilon)}$] and this concludes the proof of (i). The second assertion (ii) follows from

$$\|\Phi^{(\varepsilon)}(\mathbf{x})^{-1} - \Phi^{-1}\| = \left\|\sum_{k=1}^\infty (-\Phi^{-1}\Psi)^k\right\| \leq \frac{\bar{C}r(\varepsilon)}{1 - \bar{C}r(\varepsilon)}. \qquad \square$$

PROOF OF (4.13) IN PROPOSITION 4.5. By Proposition 4.2,

$$\frac{\partial \theta}{\partial x_i^\gamma}(\mathbf{x}) = \varepsilon^\kappa \theta(\mathbf{x}) \Phi^{(\varepsilon)}(\mathbf{x})^{-1} \operatorname{Proj}\{(\theta(\mathbf{x})^{-1}\mathfrak{e}_\gamma) \otimes z_i^{(\varepsilon)}\},$$

and therefore (4.13) follows from Lemma 5.1(i) noting $|z_i^{(\varepsilon)}| \leq R\varepsilon^{-1}$. $\square$

LEMMA 5.2. *For every* $1 \leq \gamma, \gamma' \leq d$,

$$\sup_{0<\varepsilon\leq\varepsilon_0} \sup_{\mathbf{x}\in\mathcal{M}_\infty^{(\varepsilon)}(\varepsilon^{\nu-1})} \sup_{1\leq i,i'\leq N(\varepsilon)} \varepsilon^{-2\kappa+2} \left|\frac{\partial^2 \theta}{\partial x_i^\gamma \partial x_{i'}^{\gamma'}}(\mathbf{x})\right| < \infty.$$

PROOF. The identity (4.9) may be rewritten as

$$\left(Q^{(\varepsilon)}(\mathbf{x})\frac{\partial \theta}{\partial x_i^\gamma}(\mathbf{x}) + (z_i^{(\varepsilon)} \otimes \mathfrak{e}_\gamma)\theta(\mathbf{x}), Y\right) = 0.$$

Hence, taking the derivative of this identity in $x_{i'}^{\gamma'}$, we have

$$\left(Q^{(\varepsilon)}(\mathbf{x})\frac{\partial^2 \theta}{\partial x_i^\gamma \partial x_{i'}^{\gamma'}}(\mathbf{x}) + (z_{i'}^{(\varepsilon)} \otimes \mathfrak{e}_{\gamma'})\frac{\partial \theta}{\partial x_i^\gamma}(\mathbf{x}) + (z_i^{(\varepsilon)} \otimes \mathfrak{e}_\gamma)\frac{\partial \theta}{\partial x_{i'}^{\gamma'}}(\mathbf{x}), Y\right) = 0.$$

This implies

$$\frac{\partial^2 \theta}{\partial x_i^\gamma \partial x_{i'}^{\gamma'}}(\mathbf{x}) = -\varepsilon^\kappa \theta(\mathbf{x}) \Phi^{(\varepsilon)}(\mathbf{x})^{-1}$$

$$\times \operatorname{Proj}\left\{(z_{i'}^{(\varepsilon)} \otimes \mathfrak{e}_{\gamma'})\frac{\partial \theta}{\partial x_i^\gamma}(\mathbf{x}) + (z_i^{(\varepsilon)} \otimes \mathfrak{e}_\gamma)\frac{\partial \theta}{\partial x_{i'}^{\gamma'}}(\mathbf{x})\right\}.$$

The conclusion follows from Lemma 5.1(i) and (4.13) by noting that $|z_{i'}^{(\varepsilon)}|$, $|z_i^{(\varepsilon)}| \leq R\varepsilon^{-1}$. $\square$

PROOF OF (4.14) IN PROPOSITION 4.5. Applying the mean value theorem, we have, from Lemma 5.2,

$$\left|\frac{\partial \theta}{\partial x_i^\gamma}(\mathbf{x}) - \frac{\partial \theta}{\partial x_i^\gamma}(\mathbf{z}(\mathbf{x}))\right| \leq \sum_{i'=1}^{N}\sum_{\gamma'=1}^{d} \left|\frac{\partial^2 \theta}{\partial x_i^\gamma \partial x_{i'}^{\gamma'}}(\mathbf{x}^*)\right| |x_{i'}^{\gamma'} - (\mathbf{z}(\mathbf{x}))_{i'}^{\gamma'}|$$

$$\leq CN \cdot \varepsilon^{2\kappa-2} \cdot \varepsilon^{\nu-1} \leq C'\varepsilon^{d+\nu+1},$$

where $\mathbf{x}^*$ is a certain point on the segment connecting $\mathbf{x}$ and $\mathbf{z}(\mathbf{x})$. This shows (4.14). $\square$

## 6. Concluding remarks.

6.1. We have studied in Section 4 the case where the macroscopic body is $d$-dimensional, but one can consider thin bodies and derive their motion as well. Assume that, for an $n$-dimensional Riemannian manifold $M$ in $\mathbb{R}^d, n < d$, a sequence of infinitesimally rigid crystals $\mathbf{z}^{(\varepsilon)} = (z_i^{(\varepsilon)})_{i=1}^N$ is given and it has a macroscopic limit density function $\rho_M(y)$ on $M$ in the sense that

$$(6.1) \qquad \lim_{\varepsilon \downarrow 0} \varepsilon^n \sum_{i=1}^N \delta_{\varepsilon z_i^{(\varepsilon)}}(dy) = \rho_M(y)\,dy_M,$$

where $dy_M$ denotes the volume element of $M$. Then, comparing (6.1) with the condition 2 for $\mathbf{z}^{(\varepsilon)}$ in Section 4.2, we see that $\rho(y) = 0$ for the sequence $\mathbf{z}^{(\varepsilon)}$, which we are considering here, under the scaling (4.10), and therefore $\bar{\rho} = \bar{q}^\alpha = 0$ in Theorem 4.3. This means that a different time scaling is required for $\mathbf{x}(t)$ to have a nontrivial macroscopic limit. Indeed, one can show that the right scaling is $\mathbf{x}^{(\varepsilon)}(t) := \mathbf{x}(\varepsilon^{-(n+2)}t)$, and under this scaling, random motion of the body $(M, \rho_M)$ is obtained in the limit. Note that the affine hull of $\mathbf{z}^{(\varepsilon)}$ should be $d$-dimensional to be infinitesimally rigid (see [2], page 174) and therefore, even for obtaining an $n$-dimensional macroscopic body in the limit, $d$-dimensional configurations should be considered microscopically. High polymers or membranes studied in physical chemistry usually have the above structures with $n = 1$ or $2$ in $\mathbb{R}^3$.

6.2. If the Hamiltonian is suitably modified, the notion of rigidity for the microscopic configurations may change. Let $V \in C_0^2(\mathbb{R})$ be a symmetric function having a deep well at 0 satisfying $V''(0) > 0$ and consider the modified Hamiltonian of $H(\mathbf{x})$ by adding a three-body interaction term:

$$\tilde{H}(\mathbf{x}) = H(\mathbf{x}) + \sum_{i,j,k} V\left(\left|\frac{x_i + x_j}{2} - x_k\right|\right) \mathbb{1}_{\{\mathrm{diam}\{x_i, x_j, x_k\} \leq a_0\}},$$

where the sum is taken for $i, j, k$ different and $a_0$ is a constant smaller than $4a$. We assume $b < 2a$ for the potential $U$ in $H(\mathbf{x})$. Then, one-dimensional straight chains $\mathbf{z} = (z_i)_{i=1}^N$ in $\mathbb{R}^d: z_i - z_0 = i(z_1 - z_0)$, $1 \leq i \leq N$, arranged in an equal distance $a$ (i.e., $|z_1 - z_0| = a$) are local minima of $\tilde{H}$. The proper time change is $\mathbf{x}^{(\varepsilon)}(t) := \mathbf{x}(\varepsilon^{-3}t)$, which is the same as taking $n = 1$ in Section 6.1. The particles' number $N$ should behave as $\bar{\rho} = \lim_{\varepsilon \downarrow 0} \varepsilon N$ exists. Introducing an energy different from ours, Kotani and Sunada [10] characterize the equilibrium configurations of crystals.

**Acknowledgment.** The problem studied in this paper was posed by Dmitry Ioffe. The author thanks him for valuable discussions.

# REFERENCES

[1] ASIMOW, L. and ROTH, B. (1978). The rigidity of graphs. *Trans. Amer. Math. Soc.* **245** 279–289. MR511410

[2] ASIMOW, L. and ROTH, B. (1979). The rigidity of graphs. II. *J. Math. Anal. Appl.* **68** 171–190. MR531431

[3] BODINEAU, T., IOFFE, D. and VELENIK, Y. (2000). Rigorous probabilistic analysis of equilibrium crystal shapes. *J. Math. Phys.* **41** 1033–1098. MR1757951

[4] FUNAKI, T. (1995). The scaling limit for a stochastic PDE and the separation of phases. *Probab. Theory Related Fields* **102** 221–288. MR1337253

[5] FUNAKI, T. (1997). Singular limit for reaction-diffusion equation with self-similar Gaussian noise. In *New Trends in Stochastic Analysis* (K. D. Elworthy, S. Kusuoka and I. Shigekawa, eds.) 132–152. World Scientific, Singapore. MR1654352

[6] FUNAKI, T. (2004). Zero temperature limit for interacting Brownian particles. II. Coagulation in one dimension. *Ann. Probab.* **32** 1228–1246. MR2060297

[7] FUNAKI, T. and NAGAI, H. (1993). Degenerative convergence of diffusion process toward a submanifold by strong drift. *Stochastics Stochastics Rep.* **44** 1–25. MR1276927

[8] HALES, T. C. (1997). Sphere packings. I. *Discrete Comput. Geom.* **17** 1–51. MR1418278

[9] KATZENBERGER, G. S. (1991). Solutions of a stochastic differential equation forced onto a manifold by a large drift. *Ann. Probab.* **19** 1587–1628. MR1127717

[10] KOTANI, M. and SUNADA, T. (2000). Standard realizations of crystal lattices via harmonic maps. *Trans. Amer. Math. Soc.* **353** 1–20. MR1783793

[11] KUNITA, H. (1990). *Stochastic Flows and Stochastic Differential Equations*. Cambridge Univ. Press. MR1070361

[12] LAMAN, G. (1970). On graphs and rigidity of plane skeletal structures. *J. Engrg. Math.* **4** 331–340. MR269535

[13] LANDAU, L. D. and LIFSHITZ, E. M. (1969). *Mechanics*, 2nd ed. Pergamon, New York. MR120782

[14] LANG, R. (1979). On the asymptotic behaviour of infinite gradient systems. *Comm. Math. Phys.* **65** 129–149. MR528187

[15] MCKEAN, H. P., Jr. (1969). *Stochastic Integrals*. Academic Press, New York. MR247684

[16] SHIGEKAWA, I. (1984). Transformations of the Brownian motion on a Riemannian symmetric space. *Z. Wahrsch. Verw. Gebiete* **65** 493–522. MR736143

[17] WHITELEY, W. (1984). Infinitesimally rigid polyhedra. I. Statics of frameworks. *Trans. Amer. Math. Soc.* **285** 431–465. MR752486

GRADUATE SCHOOL
OF MATHEMATICAL SCIENCES
UNIVERSITY OF TOKYO
3-8-1 KOMABA MEGURO-KU
TOKYO 153-8914
JAPAN
E-MAIL: funaki@ms.u-tokyo.ac.jp